\documentclass[11pt,epsf]{article}
\usepackage{graphicx}
\usepackage{amsthm}
\usepackage{amsfonts}
\usepackage{amssymb}
\title{Exclusion zones for the zeros of the partial theta function}
\author{Vladimir Petrov Kostov\\ 
Universit\'e C\^ote d'Azur, CNRS, LJAD, France\\   
e-mail: vladimir.kostov@univ-cotedazur.fr} 
\date{}
\newtheorem{tm}{Theorem}

\newtheorem{rem}[tm]{Remark}
\newtheorem{rems}[tm]{Remarks}
\newtheorem{lm}[tm]{Lemma}

\newtheorem{prop}[tm]{Proposition}

\begin{document}

\maketitle

\begin{abstract}
  For $q\in (-1,0)$, all complex conjugate pairs of zeros of the partial theta
  function $\theta (q,x)=\sum_{j=0}^{\infty}q^{j(j+1)/2}x^j$ belong to the
  disk $\mathbb{D}_{61.64}$. For $q\in (0,1)$, all its conjugate pairs with
  positive real parts belong to the set $\{ |x|<5$, Re$x<2.4\}$ and avoid a set containing
  the right unit half-disk. Its conjugate pairs with negative real parts are known to belong to the disk $\mathbb{D}_{49.8}$; they avoid a set
  containing the left unit half-disk and more than $12$ times larger than the latter.
  The borders of these sets consist of segments and arcs of circles.
  The function has no zeros with positive real part
  for $q\in (0,0.7234985094]$.\\


\end{abstract}

\section{Introduction}

\subsection{Formulation of the results}

The partial theta function is the sum of the series
$\theta (q,x)=\sum_{j=0}^{\infty}q^{j(j+1)/2}x^j$ in which $x$ is a variable and $q$
a parameter. We are interested in the case when $x$ is real and
$q\in (-1,0)\cup (0,1)$. The applications of this function are briefly
discussed in the next subsection.

The name of the partial theta function is derived from the one of the
{\em Jacobi theta function} $\Theta (q,x):=\sum_{j=-\infty}^{\infty}q^{j^2}x^j$,
having in mind the equality $\theta (q^2,x/q)=\sum_{j=0}^{\infty}q^{j^2}x^j$. ``Partial''
refers to the fact that in the case of $\theta$ summation is performed only
on the non-negative integers.

The Jacobi theta function can be represented
via the {\em Jacobi triple product}:

\begin{equation}\label{equtriple}
  \Theta (q,x^2)=\prod_{m=1}^{\infty}(1-q^{2m})(1+x^2q^{2m-1})(1+x^{-2}q^{2m-1})~.
  \end{equation}
When applied to the function

$$\Theta^*(q,x)=\Theta (\sqrt{q},\sqrt{q}x)=
\sum_{j=-\infty}^{\infty}q^{j(j+1)/2}x^j~,$$
formula (\ref{equtriple}) gives its representation

\begin{equation}\label{equT}
  \begin{array}{ccl}
    \Theta^*(q,x)&=&\prod_{m=1}^{\infty}(1-q^m)(1+xq^m)(1+q^{m-1}/x)\\ \\
    &=&(1+1/x)\prod_{m=1}^{\infty}(1-q^m)(1+xq^m)(1+q^m/x)~.\end{array}
\end{equation}
We define the function
$G:=\sum_{j=-\infty}^{-1}q^{j(j+1)/2}x^j=1/x+q/x^2+q^3/x^3+q^6/x^4+\cdots$ to write

\begin{equation}\label{equthetaG}
  \theta =\Theta^*-G~.
  \end{equation}
The partial theta function satisfies the functional equation

\begin{equation}\label{equFE}
  \theta (q,x)=1+qx\theta (q,qx)~.
  \end{equation}
For each $q\in (0,1)$ fixed, the function $\theta$ has at most finitely-many 
complex conjugate pairs, infinitely-many negative and no positive zeros. For
each $q\in (-1,0)$, the partial theta function has finitely-many
conjugate pairs and infinitely-many
positive and negative zeros. Real zeros go to infinity asymptotically as
$-1/q^n$.
See more about the analytic properties of $\theta$ in
Subsection~\ref{subsecmore}. At the same time (see Theorem~\ref{tmoldresults})
complex conjugate
pairs remain within a compact set independent of~$q$.
Its complement in $\mathbb{C}$ could be termed as {\em exclusion zone}. 
The aim of the present paper is to prove new and better results about the  
location of the zeros which is tantamount to expanding the exclusion zone.
We remind 
what is known in this sense:

\begin{tm}\label{tmoldresults}
(1) {\rm (\cite[Theorem~1]{KoArxiv})}
  For $q\in (0,1)$,
  the complex conjugate pairs with non-negative real part (if any) of
  $\theta (q,.)$ belong to the half-annulus
  $\mathcal{A}:=\{ {\rm Re}x\geq 0,~1<|x|<5\}$.
  \vspace{1mm}
  
  (2) {\rm (\cite[Theorem~3]{KoArxiv})} For $q\in (0,1)$, the complex
  conjugate pairs
  of zeros of $\theta (q,.)$ with negative real part
belong to the left open half-disk of radius $49.8$ centered at the origin.
\vspace{1mm}

(3) {\rm (\cite[Theorem~1]{KoMatStud})} For any fixed $q\in (0,1)$,
the partial theta function has no zeros in the domain
$\mathcal{D}:=\{ x\in \mathbb{C}:|x|\leq 3,~{\rm Re} x\leq 0,
|{\rm Im} x|\leq 3/\sqrt{2}\}$ (with $3/\sqrt{2}=2.121320344\ldots$).
\vspace{1mm}

(4) {\rm (\cite[Theorem~1]{KoAnn24})} For each $q\in (-1,0)\cup (0,1)$ fixed,
the function $\theta$ has no zeros in the closed unit disk
$\overline{\mathbb{D}_1}$.
\vspace{1mm}

(5) {\rm (\cite[Theorem~1]{KoFAA19})} For $q\in (-1,0)$, all complex conjugate
zeros of the partial theta function belong to the set
$\{ x\in \mathbb{C},~|{\rm Re}x|\leq 364.2,~|{\rm Im}x|\leq 132\}$. 
\end{tm}

\begin{figure}[htbp]
\centerline{\hbox{\includegraphics[scale=0.7]{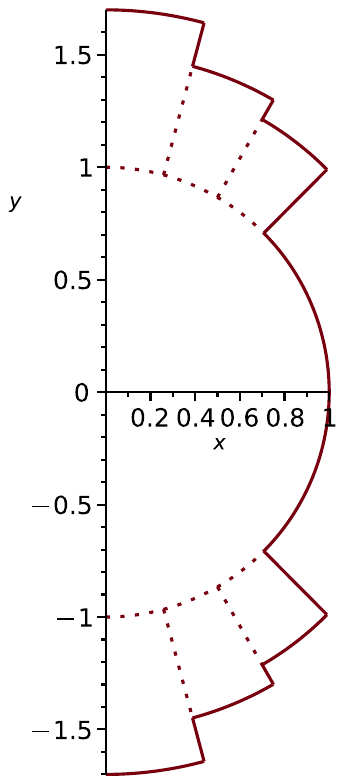}}}
\caption{The contour of the domain $\tilde{R}$ (without the segment $S_1$) in
  solid line and the rest of the right unit circle and the continuations of
  the segments $S_2$, $S_3$, $S_4$ and their symmetric ones w.r.t.
  the abscissa-axis in dotted line.}
\label{figfresnel1}
\end{figure}

\begin{figure}[htbp]
\centerline{\hbox{\includegraphics[scale=0.7]{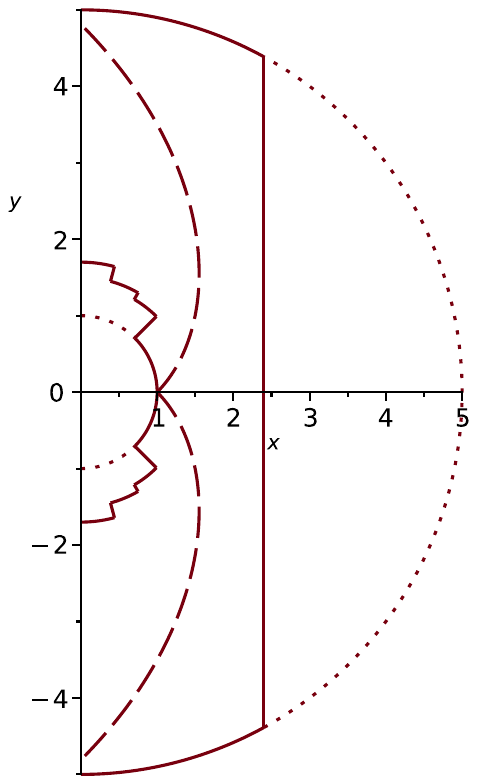}}}
\caption{The border of the set $\mathcal{B}$ (see Remarks~\ref{remsestim})
    in solid line, the rest of the unit circle
  and the circle centered at $0$ and of radius $5$ in dotted line and the
intersection of Katsnelson's contour with the right half-plane in dashed line.}
\label{figfresnel2}
\end{figure}

Our first new result reads:

\begin{tm}\label{tmlocation}
  For $q\in (-1,0)$, all complex conjugate pairs of zeros of $\theta$
  belong to the open disk $\mathbb{D}_{61.64}$. 
\end{tm}
The theorem is proved in Section~\ref{secprtmlocation}. The area of the disk
$\mathbb{D}_{61.64}$ is more than $16$ times smaller than the one of the
rectangle mentioned in part (5) of Theorem~\ref{tmoldresults}.
In order to formulate
the next result we define the closed domain $\tilde{R}\subset \mathbb{C}$. It is
bounded by the arcs

$$\begin{array}{cclccl}
  C_1&:=&\{ |x|=1.7,~{\rm arg}x\in [5\pi /12,\pi /2]\}~,&
  C_2&:=&\{ |x|=1.5,~{\rm arg}x\in [\pi /3,5\pi /12]\}~,\\ \\ 
  C_3&:=&\{ |x|=1.4,~{\rm arg}x\in [\pi /4,\pi /3]\}~,&
  C_4&:=&\{ |x|=1,~{\rm arg}x\in [0,\pi /4]\}~,\end{array}$$
their symmetric ones w.r.t. the abscissa-axis, the segment $S_1:=[-1.7i,1.7i]$
and the segments

$$S_2:=[1.5e^{5\pi i/12},1.7e^{5\pi i/12}]~,~~~\,
S_3:=[1.4e^{\pi i/3},1.5e^{\pi i/3}]~,~~~\, S_4:=[e^{\pi i/4},1.4e^{\pi i/4}]$$
and their symmetric ones w.r.t. the abscissa-axis, see Fig.~\ref{figfresnel1}
and Fig.~\ref{figfresnel2}.

\begin{figure}[htbp]
\centerline{\hbox{\includegraphics[scale=0.7]{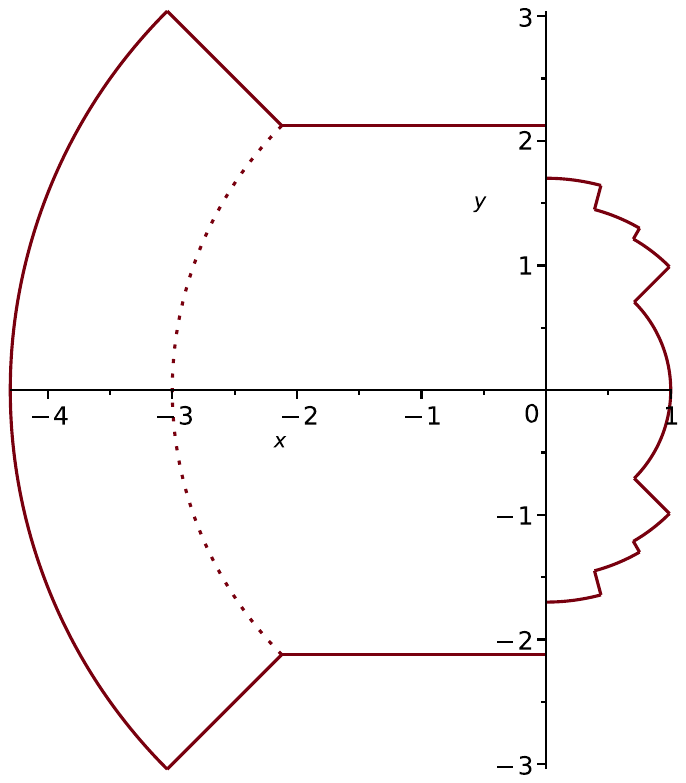}}}
\caption{The exclusion zones $\tilde{R}$, $\mathcal{D}$ and $\mathcal{F}$ around the origin.}
\label{figleftrightexclusion}
\end{figure}

\begin{figure}[htbp]
\centerline{\hbox{\includegraphics[scale=0.7]{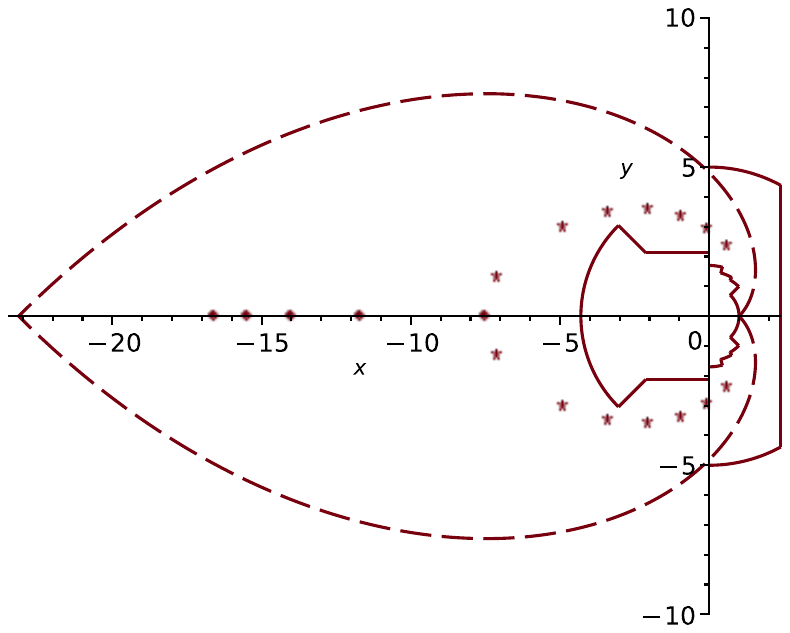}}}
\caption{The exclusion zones around the origin and in the right half-plane (in solid line), Katsnelson's contour (in dashed line), the first five double zeros ($-7.5\ldots$, $-11.7\ldots$, $-14.0\ldots$, $-15.5\ldots$ and $-16.6\ldots$, corresponding to the spectral values $0.309249\ldots$, $0.516959\ldots$, $0.630628\ldots$, $0.701265\ldots$ and $0.749269\ldots$) and the complex conjugate pair (computed for an approximation of $\theta$ with 61 terms) born from the double zero $-7.5\ldots$ for $q=0.32$, $0.4$, $0.48$, $0.56$, $0.64$, $0.72$ and $0.8$ (in asterisks).}
\label{figexclusion}
\end{figure}

\begin{tm}\label{tmimprove}
  (1) For $q\in (0,1)$, the complex
  conjugate pairs of zeros of $\theta (q,.)$ with positive real part
  belong to the half-plane {\rm Re}$x<2.4$.
\vspace{1mm}

  (2) For $q\in (0,1)$, the complex
  conjugate pairs of zeros of $\theta (q,.)$ with non-negative real part
  avoid the set $\tilde{R}$. 

\end{tm}

The theorem is proved in Section~\ref{sectmimprove}.

\begin{tm}\label{tmimprovebis}
 For $q\in (0,1)$, there are no zeros of $\theta$ in the quarter of annulus $\mathcal{F}:=\{ 3\leq |x|\leq 4.3$, {\rm arg}$x\in [3\pi /4,5\pi /4]\}$.
\end{tm}

The theorem is proved in Section~\ref{sectmimprovebis}. The sets $\tilde{R}$, $\mathcal{D}$ and $\mathcal{F}$ are shown in Fig.~\ref{figleftrightexclusion}; the common border of $\mathcal{D}$ and $\mathcal{F}$ is drawn in dotted line. The area of the set $\mathcal{D}\cup \mathcal{F}$ is $19.0220\ldots$ which is 
$12.109\ldots$ times more than the area of the left unit half-disk.

\begin{rems}\label{remsestim}
  {\rm (1) The area of the domain
    $\tilde{R}$ is $68.3\%$ larger than the area of the right unit half-disk
    (see about the half-annulus $\mathcal{A}$ in part (1) of
    Theorem~\ref{tmoldresults}). These are the exclusion zones in the right
    half-plane around the origin for the complex conjugate zeros of $\theta$
    according to Theorems~\ref{tmoldresults} and \ref{tmimprove}. 
    \vspace{1mm}
    
    (2) The domain in the right half-plane to which according to
    Theorem~\ref{tmimprove} the complex conjugate pairs with non-negative real
    parts belong, is} 
    
  $$\mathcal{B}:=\{ |x|<5,~{\rm Re}x\geq 0\}\setminus
    (\tilde{R}\cup
      \{ {\rm Re}x\geq 2.4\} )~.$$
    {\rm Its area is 54,11\% of the area of the half-annulus
    $\mathcal{A}$.
    \vspace{1mm}
    
    (3) In \cite{Ka} V.~Katsnelson has proved that as $q\rightarrow 1^-$, the
    series of $\theta$ converges to $1/(1-x)$ inside {\em Katsnelson's contour}
    $\{ {\rm Re}x=e^t\cos t$, Im$x=\pm e^t\sin t $, $t\in [0,\pi ]\}$. The
    leftmost point of the contour is $-e^{\pi}=-23.14\ldots$, see
    Fig.~\ref{figexclusion}. The exclusion zone around $0$, like the interior
    of Katsnelson's contour, is larger in the left and narrower in the right
    half-plane.
    \vspace{1mm}
    
    (4) In Fig.~\ref{figexclusion} one can see the complex conjugate pair of
    the approximation $\sum_{j=0}^{60}q^{j(j+1)/2}x^j$ of $\theta$ which is close to
    the pair of $\theta$ born from the double zero of $\theta (\tilde{q}_1,.)$.
    This pair is shown for $q=0.32+k\cdot 0.08$, $k=0$, $\ldots$, $6$. For
    $q=0.48$, $0.72$ and $0.8$, it is close to the exclusion zone around $0$.
    If one finds in a similar way the positions of the pairs born from the
    double zeros for $q=\tilde{q}_j$, $j=2$, $3$ and $4$, one will see that the
    pairs go to the right  each staying closer  to Katsnelson's contour than
    the previous one. Thus one can suggest that complex conjugate zeros remain
    in the interior of the contour and that if one restrains $q$ to an interval
    of the form $[a,1)$, $0<a<1$, then one can find an exclusion zone around
      $0$ which expands as $a$ increases and which  tends to the interior of
      Katsnelson's contour as $a\rightarrow 1^-$. One can expect that as $q$
      grows, each complex pair enters (for some value of $q$, individual for
      each pair) the right half-plane and then tends to $1$.}
\end{rems}

The following proposition
contains a functional equation satisfied by $\theta$; it is
to be checked straightforwardly:

\begin{prop}\label{propdecomp}

  For $n\in \mathbb{N}^*$, one has
  $\theta (q,x)=\sum_{j=0}^{n-1}q^{j(j+1)/2}x^j\theta (q^{n^2},q^{-n(n-1)/2+nj}x^n)$.
\end{prop}

We use in the proofs of the theorems also a proposition
(proved in Section~\ref{secpropheart}) which is of
independent interest:

\begin{prop}\label{propheart}
  For $q\in (0,q_{\heartsuit}]$, $q_{\heartsuit}:=0.7234985094$, and
      ${\rm Re}x\geq 0$, one has
    $\theta (q,x)\neq 0$.
\end{prop}

\subsection{Applications of the partial theta function}

The partial theta function has found applications in the following domains:

\begin{itemize}
\item
  modularity and asymptotics of 
regularized characters and 
quantum dimensions of the $(1,p)$-singlet algebra modules (see \cite{BFM} 
and \cite{CMW});

\item
  the theory 
  of (mock) modular forms (see \cite{BrFoRh});

\item
  Ramanujan-type $q$-series 
  (see \cite{Wa});

\item
  asymptotic analysis (see \cite{BeKi});

\item
  statistical physics 
  and combinatorics (see \cite{So});

\item
  thermodynamic phases in first detected return times
  of quantum many-body systems, see \cite{WPG}.
  \end{itemize}
Its different properties have been investigated in \cite{AnBe}, \cite{EM},
\cite{FGM}, \cite{LuSa}, \cite{Pr}, \cite{Sun}, \cite{WM}, \cite{Wei}
and the list is not exhaustive. 

Classical results of Hardy, Petrovitch and 
Hutchinson (see \cite{Ha}, \cite{Pe} and \cite{Hu}) about
{\em section-hyperbolic 
polynomials}, i.e. real univariate polynomials with all 
roots real and such that all their truncations of positive degree have also 
all their roots real, have given an impetus to a more recent interest in the
partial theta function,
see \cite{Ost}, \cite{KaLoVi} and \cite{KoSh}.

The analytic properties of the function $\theta$ have been studied in
\cite{KoFAA}, \cite{KoDBAN14}, \cite{KoDBAN15}, \cite{KoBSM1}, \cite{KoBSM2},
\cite{KoPMD}, \cite{KoAM}, \cite{KoSe} and~\cite{Kosepar}.

\subsection{Analytic properties of the partial theta function
  \protect\label{subsecmore}}

For $q\in (0,1)$, there exists a sequence of countably-many
{\em spectral} values of $q$ denoted by
$0<\tilde{q}_1<\tilde{q}_2<\cdots$, where $\tilde{q}_1=0.3092\ldots$ and
$\lim_{j\rightarrow \infty}\tilde{q}_j=1^-$. They enjoy the following properties:

\begin{itemize}
\item
  For $q=\tilde{q}_j$, the function
  $\theta$ has exactly one double negative zero which is the rightmost of
  the negative zeros; its other negative zeros
  are simple.

\item
  For the other values of $q$, the function $\theta$ has no
  multiple real zeros.

\item
  For $q\in (\tilde{q}_j,\tilde{q}_{j+1}]$, the function $\theta$ has exactly
    $j$ complex conjugate pairs of zeros counted with multiplicity;
    we set $\tilde{q}_0:=0$.
\end{itemize}

For $q\in (-1,0)$, the partial theta function has also spectral values.
They are denoted by $-1<\bar{q}_{j+1}<\bar{q}_j<0$, see the first eight of them
in (\ref{equtable}).

\begin{itemize}
\item
  For such values, $\theta$ has one double, 
infinitely-many positive and infinitely-many negative simple real zeros.

\item
  For odd $j$, the double zero is the rightmost of the negative zeros and
  $\theta$ has a local minimum there.

\item
  For even $j$, the double zero is the second from the left of the
  positive zeros and $\theta$ has a local maximum there.
  These properties of the spectrum are proved in
  \cite{KoPRSE2} and \cite{Kosepar}.

\item
  For $q\in [\bar{q}_{j+1},\bar{q}_j)$ (with $\bar{q}_0=0$), $\theta$ has
    $j$ complex conjugate pairs of zeros counted with multiplicity. Complex
    zeros do not cross the
    imaginary axis as $q$ decreases from $0$ to $-1$,
    see~\cite[Theorem~8]{KoSe}.
  \end{itemize}

\section{Proof of Proposition~\protect\ref{propheart}\protect
  \label{secpropheart}}

The zeros of $\theta$ depend continuously on the parameter $q$.
For $q\in (0,\tilde{q}_1]$, all zeros are negative. Thus it suffices to show 
that for $q\in (0,q_{\heartsuit}]$, no zero of $\theta$ crosses the imaginary 
axis (this cannot happen at $0$, because $\theta (q,0)=1$). For $x=iy$, 
$y\in \mathbb{R}$, we use the identity 

      $$\theta (q,iy)=v_R+iqyv_I~,~~~\,
      v_R:={\rm Re}\theta =\theta (q^4,-y^2/q)~,~~~\,
      v_I:={\rm Im}\theta =\theta (q^4,-qy^2)~.$$
      We have to prove that for $q\in (0,q_{\heartsuit}]$, the functions $v_R$ and
        $v_I$ have no zero in common. We set $w:=q^4$ and $\xi :=-y^2/q$. Thus
        $v_R=\theta (w,\xi )$ and $v_I=\theta (w,\sqrt{w}\xi )$.  
        
   Consider the functions $\tau_1(w):=\theta (w,-w^{-1.4})$ and
  $\tau_2(w):=\theta (w,-w^{-1.9})$. One can show numerically that these two
  functions are both negative for $w\in (0,0.274]$. We apply to the function
    $\theta (q,x)$ formula 
    (\ref{equFE}) twice:

    $$\theta (q,x)=1+qx+q^3x^2\theta (q,q^2x)~.$$
    Thus for $k\in 2\mathbb{N}$,

    $$\begin{array}{ccl}
      \theta (w,-w^{-1.4-k})&=&1-w^{-0.4-k}+w^{0.2-2k}\theta (w,-w^{0.6-k})~~~\,
      {\rm and}\\ \\ 
      \theta (w,-w^{-1.9-k})&=&1-w^{-0.9-k}+w^{-0.8-2k}\theta (w,-w^{0.1-k})~.
      \end{array}$$
        The above equalities allow to prove by induction on $k$ that 
        $${\rm if}~~~\, \theta (w,-w^{-\beta })<0~,~~~\, 
        {\rm then~~~\, }\theta (w,-w^{-\beta -k})<0~,~~~\, \beta =1.4~~{\rm or}~~
        1.9~.$$
        Indeed, if $\theta (w,-w^{0.6-k})<0$, then as $w^{-0.4-k}>1$, one obtains
        $\theta (w,-w^{-1.4-k})<0$ (and one  similarly proves that
        $\theta (w,-w^{-1.9-k})<0$). Knowing that $\theta (w,-w^{-n})>0$,
        $n\in \mathbb{N}^*$ (\cite[Proposition~9]{KoBSM1}),
        $\theta (w,-w^{-2s-1.4})<0$ and $\theta (w,-w^{-2s-1.9})<0$,
        $s\in \mathbb{N}$,  one concludes that
        the zeros of the function $v_R=\theta (w,\xi )$ belong to the intervals 
        $$(-w^{-2s-1.4},-w^{-2s-1})~~~\, {\rm 
        and}~~~\, (-w^{-2s-2},-w^{-2s-1.9})~,~~~\, s\in \mathbb{N}~.$$
        Then 
        the zeros of the function $v_I=\theta (w,\sqrt{w}\xi )$ belong to the intervals
        
        $$(-w^{-2s-1.9},-w^{-2s-1.5})~~~\, {\rm and}~~~\,
        (-w^{-2s-2.5},-w^{-2s-2.4})~.$$
        These two
        couples of families of intervals do not intersect. Therefore for
        $w\in (0,0.274]$
          (hence for $q\in (0,0.274^{1/4}]$,
            $0.274^{1/4}=0.7234985094\ldots >q_{\heartsuit}$),
          the functions
          $v_R$ and $v_I$ have no zero in common, so $\theta \neq 0$.

      \section{Proof of Theorem~\protect\ref{tmlocation}
            \protect\label{secprtmlocation}}

The idea of proof looks as follows. We show that $\theta (q,x)\neq 0$
along certain circles in the plane of the complex variable $x$
centered at $0$. Their radii decrease as $q$ decreases
in $(-1,0)$. The complex conjugate pair born from the double zero of
$\theta (\bar{q}_1,.)$ belongs to the interior of the first circle. For
$q=\bar{q}_2$, we change the first circle for the second one whose radius is
larger and whose interior contains the first two complex conjugate pairs. Then
for $q=\bar{q}_3$, we change the second circle for the third one whose
interior contains the first three complex conjugate pairs etc. The radius of the
circle containing in its interior all conjugate pairs
is a piecewise-continuous function in $q$ whose values are $<61.64$. Changing a circle by a larger one is compensated by the decreasing of the radius between two changes, so the radius is never larger than $61.64$. 
The idea is developed in Subsections~\ref{subsecfirstseveral}
and~\ref{subsecnextspectral}. In the next two
subsections we present the necessary technical tools.

\subsection{Beginning of the proof}

For $q\in (-1,0)$ and Re$x\leq 0$, consider the product $\Theta^*(q,x)$. Its
factors with odd powers of $q$ define the quantity

$$\Theta ^{**}(q,x)=\prod_{m=1}^{\infty}\left( (1-q^{2m-1})(1+q^{2m-1}x)
(1+q^{2m-1}/x)\right) ~.$$
It is clear that

$$1-q^{2m-1}>1~,~~~\, {\rm Re}(1+q^{2m-1}x)\geq 1~~~\, {\rm and}~~~\, 
{\rm Re}(1+q^{2m-1}/x)\geq 1~,$$
therefore $|\Theta^{**}|>1$. It is also evident that

\begin{equation}\label{equT*}
  \Theta^*(q,x)=\Theta^{**}(q,x)\Theta^*(q^2,x)~,~~~\, {\rm so}~~~\,
  |\Theta^*(q,x)|>|\Theta^*(q^2,x)|~.
  \end{equation}
In what follows we write $q_*$ instead of $q$ when considering
case A)~:~$q\in (0,1)$;
we write $q$ when we are in case B)~:~$q\in (-1,0)$. 
We observe that for both cases, A) and B), and for $|x|\geq 5$, one obtains
the estimation
\begin{equation}\label{equ|G|}
  |G|\leq \sum_{j=1}^{\infty}|q^{j(j-1)/2}|/|x|^j<
  \sum_{j=1}^{\infty}|x|^{-j}\leq \sum_{j=1}^{\infty}5^{-j}=1/4~.
  \end{equation}
The lemma that follows is a more precise version of Lemma~16 in~\cite{KoArxiv}:

\begin{lm}\label{lm-5}
  Suppose that for $q_*\in (0,1)$ fixed and for $x=x_0<-5$, one has
  $|\theta (q_*,x_0)|\geq 1$. Then one has $|\theta (q_*,z)|>11/20$ for
  $|z|=|x_0|$.
  \end{lm}

\begin{proof}[Proof of Lemma~\ref{lm-5}]
  For $z=x_0$, as $x_0<0$, the series $G$ is a Leibniz series hence the modulus
  of its sum is majorized by the modulus of its first term,
  so $|G|<1/5$ and $0<-G<1/5$.

  In the proof of Lemma~16 in \cite{KoArxiv} it is shown that for $|z|=|x_0|$,
  the modulus $|\Theta^*(q_*,z)|$ is minimal for $z=x_0$. Recall that
  $\theta =\Theta^*-G$. Suppose that $\theta (q_*,x_0)\leq -1$. Then

  $$\Theta^*(q_*,x_0)=\theta (q_*,x_0)+G(q_*,x_0)<-1~~~\, {\rm and}~~~\,
  |\Theta^*(q_*,z)|>1~~~\, {\rm for}~~~\, |z|=|x_0|~.$$
  From $|G(q_*,z)|<1/4$ (see (\ref{equ|G|})) one concludes that
  $|\theta (q_*,z)|>3/4>11/20$.

  Suppose now that
  $\theta (q_*,x_0)\geq 1$. Then $\Theta^*(q_*,x_0)>4/5$, so
  $|\Theta^*(q_*,z)|>4/5$ and $|\theta (q_*,z)|>4/5-1/4=11/20$ for $|z|=|x_0|$.
  \end{proof}

We explain how we apply Lemma~\ref{lm-5} (with $q_*=q^2$) in the proof of
Theorem~\ref{tmlocation}. From equality (\ref{equthetaG}), the inequality 
$|\theta (q^2,z)|>11/20$ with $|z|=|x_0|$, 
the inequality (\ref{equT*}) and the estimation 
(\ref{equ|G|}) one deduces that

\begin{equation}\label{equestimation}
  \begin{array}{ccccl}|\theta (q,z)|&\geq&|\Theta^*(q,z)|-|G(q,z)|
  &>&|\Theta^*(q^2,z)|-|G(q,z)|\\ \\ 
    &>&|\theta (q^2,z)|-|G(q^2,z)|-|G(q,z)|&>&11/20-1/4-1/4=1/20>0~.\end{array}
\end{equation}
Hence there is no zero of $\theta$ with $|z|=|x_0|$.

\subsection{Inequalities involving the zeros of $\theta$}

We denote by $x_j$ the zeros of $\theta (q,.)$ for $q\in (-1,0)$
and we remind that (see~\cite[Fig.~3]{KoPRSE2}) whenever the indicated zeros
are real and distinct, they satisfy the following string of inequalities:

\begin{equation}\label{equineq}
  \begin{array}{l}\cdots <x_{4k+6}<qx_{4k+7}
    <qx_{4k+5}<x_{4k+4}<x_{4k+2}<qx_{4k+3}<\cdots <0\\ \\
    0<\cdots <x_{4k+1}<qx_{4k+2}<qx_{4k+4}<x_{4k+3}<x_{4k+5}<qx_{4k+6}<\cdots ~.
\end{array}\end{equation}
The first several
  spectral values $\bar{q}_j\in (-1,0)$, the
  corresponding double zeros of $\theta$ (see~\cite{KoPRSE2}) and the
  zeros that coalesce are given in the table below:

  \begin{equation}\label{equtable}
    \begin{array}{llllllll}
  -0.72\ldots&-0.78\ldots&-0.84\ldots&-0.86\ldots&-0.88\ldots&
  -0.89\ldots&-0.91\ldots&-0.91\ldots \\ \\
  -2.991\ldots&2.907\ldots&-3.621\ldots&3.523\ldots&-3.908\ldots&
  3.823\ldots&-4.08\ldots&4.002\ldots\\ \\
  x_2,x_4&x_3,x_5&x_6,x_8&x_7,x_9&x_{10},x_{12}&x_{11},x_{13}&x_{14},x_{16}&
  x_{15},x_{17}\end{array}
  \end{equation}

  \begin{rem}\label{remsigns}
    {\rm One has $\theta (q,x)>0$ for $x\in (x_{4k+2},x_{4k})$ and for
  $x\in (x_{4k+3},x_{4k+5})$,
      and $\theta (q,x)<0$ for $x\in (x_{4k+4},x_{4k+2})$ and for
      $x\in (x_{4k+1},x_{4k+3})$.}
    \end{rem}

\begin{lm}\label{lm2zeros}
  (1) For $q\in (-1,0)$ and when the zeros $x_{4k+2}$ and $x_{4k+4}$ are real,
  one has $-q^{-4k-4}<x_{4k+4}<x_{4k+2}<-q^{-4k-2}$ and 
  $\theta (q,-q^{-2\ell})>0$, $k$, $\ell \in \mathbb{N}$.

  (2) For $q\in (-1,0)$ and $\ell \in \mathbb{N}$, one has
  $|\theta (q,-q^{-2\ell -1})|<1$.
\end{lm}

\begin{proof}
  Part (1). For $\nu \in \mathbb{N}$, it is true that

  \begin{equation}\label{equthetanu}
    \theta (q,-q^{-\nu})=\sum_{j=0}^{\infty}(-1)^jq^{-\nu j+j(j+1)/2}=
    \sum_{k=1}^{\infty}(-1)^{k-1}q^{\nu k+k(k-1)/2}=:S_{\dagger}
    \end{equation}
  (to be checked directly; in the series of $\theta (q,-q^{-\nu})$
  the first term cancels with the $(2\nu )$-th,
  the second with the $(2\nu -1)$-st etc.). The change of the sommation index
  $k=j-2\nu +1$ gives then formula~(\ref{equthetanu}).

  If $\nu =2\ell$,
  $\ell \in \mathbb{N}$, then the even and odd terms of the series
  $S_{\dagger}$ form two Leibniz series $S_1$ and $S_2$ of the form
  $\pm (|q|^{a_1}-|q|^{a_2}+|q|^{a_3}-\cdots )$, where $a_j\in \mathbb{N}^*$,
  $a_1<a_2<a_3<\cdots$; in the context of part (1) of the lemma 
  their first terms are positive. This
  proves the second claim of part (1) of the lemma. For $q\in (\bar{q}_1,0)$, all zeros
  $x_j$ are real. For $q\rightarrow 0^-$, the zero $x_{4k+2}$ (resp. $x_{4k+4}$)
  is close to $-q^{-4k-2}$ (resp. $-q^{-4k-4}$); this follows from the
  asymptotic expansions of the zeros of $\theta$ as Laurent series in $q$,
  see~\cite{KoDBAN15}.
  This proves the first claim
  of part~(1).
  \vspace{1mm}
   
  Part (2). If $\nu =2\ell +1$, then the first terms of the two Leibniz series
  $S_1$ and $S_2$ have opposite signs. Hence
  $|S_1|<1$, $|S_2|<1$ and as sgn$(S_1)=-$sgn$(S_2)$, one has $|S_1-S_2|<1$.
  \end{proof}

\subsection{The first several spectral values
  \protect\label{subsecfirstseveral}}

$1^{\circ}$. We construct a piecewise-continuous
    function $x_{\natural}$ with the property the open disk centered at $0$ and
    of radius $|x_{\natural}(q)|$ to contain all complex conjugate pairs of zeros
    of $\theta (q,.)$. The jump points of $|x_{\natural}|$ are at $\bar{q}_j$,
    $j=2$, $\ldots$, $6$; the function $|x_{\natural}|$ is continuous from the
    left and at a jump point its left limit is larger than its right limit.
We explain in detail the construction for $q\in [\bar{q}_3,\bar{q}_1]$ (for
    $q\in [\bar{q}_1,0)$, there are no complex conjugate pairs).
  \vspace{1mm}
  
$2^{\circ}$. Suppose that $q\in [\bar{q}_2,\bar{q}_1]$. We show that one can find
$x_0\in (-61.64,-5)$ with $|\theta (q,x_0)|>1$. We use the equality
(\ref{equFE}). If $x=x_{\natural ,1}\in (-q^{-10},-q^{-8})$, then

  $$\begin{array}{l}
    qx_{\natural ,1}\in (-q^{-7},-q^{-9})\subset (|(\bar{q}_2)^{-7}|,
  |(\bar{q}_1)^{-9}|)
  =(5.50\ldots ,17.59\ldots )~~~\, {\rm and}\\ \\ 
  x_{\natural ,1}\in (-24.20\ldots ,-7.02\ldots )~.\end{array}$$
  We choose $x_{\natural ,1}$ also such that $qx_{\natural ,1}\in (x_7,x_9)$ which
  is possible; indeed, for $q\in [\bar{q}_2,\bar{q}_1]$, the zeros $x_7$
  and $x_9$ are real and distinct and one
  has $qx_8<x_7<x_9<qx_{10}$, see (\ref{equineq}) with $k=1$. Hence
  $\theta (q,qx_{\natural ,1})>0$.
  In addition we choose $x_{\natural ,1}$ to depend continuously on
  $q\in [\bar{q}_2,\bar{q}_1]$. Thus
  (see (\ref{equFE})) $\theta (q,x_{\natural ,1})>1$ and one has
  $|\theta (q,z)|>11/20$ for $|z|=|x_{\natural ,1}|$ (Lemma~\ref{lm-5}).

  For $q=\bar{q}_1$, the
  double zero of $\theta$ equals $-2.991\ldots$, see~\cite{KoPRSE2}. Thus for
  $q\in [\bar{q}_2,\bar{q}_1]$, the only complex conjugate pair of
    $\theta (q,.)$ remains within the disk
  $\mathbb{D}_{|x_{\natural ,1}|}\subset \mathbb{D}_{61.64}$, so we can set
  $x_0=x_{\natural ,1}$.
\vspace{1mm}

$3^{\circ}$. Suppose that $q\in [\bar{q}_3,\bar{q}_2]$. The double zero of
$\theta (\bar{q}_2,.)$ equals $2.907\ldots$, see table~(\ref{equtable}).
We choose

  $$x_{\natural ,2}\in (-q^{-14},-q^{-12})\subset
  (-(\bar{q}_2)^{-14},-(\bar{q}_3)^{-12})=(-30.30\ldots ,-7.91\ldots )$$
  depending continuously on $q\in [\bar{q}_3,\bar{q}_2]$ and
  such that $qx_{\natural ,2}\in (x_{11},x_{13})$, so
  $\theta (q,qx_{\natural ,2})>0$. Hence $\theta (q,x_{\natural ,2})>1$
  (see (\ref{equFE})),
  $|\theta (q,z)|>11/20$ for $|z|=|x_{\natural ,2}|$ (see Lemma~\ref{lm-5}),
  therefore the complex conjugate pair of $\theta$ born for $q=\bar{q}_2$
  is inside the
  disk $\mathbb{D}_{|x_{\natural ,2}|}\subset \mathbb{D}_{61.64}$
  and we can set $x_0:=x_{\natural ,2}$.
For $q=\bar{q}_2$, one has

$$x_{\natural ,1}\in (-(\bar{q}_2)^{-10},-(\bar{q}_2)^{-8})~~~\, 
{\rm and}~~~\, 
x_{\natural ,2}\in (-(\bar{q}_2)^{-14},-(\bar{q}_2)^{-12})~,$$
so $|x_{\natural ,2}|>|x_{\natural ,1}|$. Hence the first complex conjugate pair of
  $\theta$ also belongs to the disk
  $\mathbb{D}_{|x_{\natural ,2}|}\subset \mathbb{D}_{61.64}$.
\vspace{1mm}

$4^{\circ}$. Now we explain the construction of $x_{\natural}$
for $q\in (\bar{q}_7,\bar{q}_1]$; for $q\in (-1,\bar{q}_7]$, the construction 
is explained in the next subsection. We constructed continuous functions
$x_{\natural ,1}$ and $x_{\natural ,2}$ on the intervals $[\bar{q}_2,\bar{q}_1]$ and
$[\bar{q}_3,\bar{q}_2]$
which share the endpoint $\bar{q}_2$. In the same way the
construction of the functions $x_{\natural ,k}$
can be carried out for $q\in [\bar{q}_{k+1},\bar{q}_k]$,
$k=3$, $4$, $5$ and $6$. Then we set $x_{\natural}:=x_{\natural ,k}$ for
  $q\in (\bar{q}_{k+1},\bar{q}_k]$, $k=1$, $\ldots$, $6$. Thus 
  the function $x_{\natural}$ is piecewise-continuous. Its jump points are at
  $\bar{q}_j$, $j=2$, $\ldots$, $6$, where it is continuous from the left;
  at a jump point its left limit is larger than its right limit. We give
  below the basic properties of the construction.

  \begin{enumerate}
  \item
    We require that $x_{\natural ,k}\in (-q^{-4k-6},-q^{-4k-4})$, so

    $$\begin{array}{ccl}
      qx_{\natural ,k}&\in&(-q^{-4k-3},-q^{-4k-5})\subset (|(\bar{q}_{k+1})^{-4k-3}|,
    |(\bar{q}_k)^{-4k-5}|)=:(\alpha_k,\beta_k)~~~\, {\rm and}\\ \\ 
    x_{\natural ,k}&\in &(-|(\bar{q}_k)^{-4k-6}|,-|(\bar{q}_{k+1})^{-4k-4}|)=:
    (\gamma_k,\delta_k)~.\end{array}$$
    We list the 2-digit truncations of the
    intervals $(\alpha_k,\beta_k)$ and $(\gamma_k,\delta_k)$ (the
    first line is $(\alpha_1,\beta_1)$, $(\gamma_1,\delta_1)$,
    $(\alpha_2,\beta_2)$, $(\gamma_2,\delta_2)$):
    
    \begin{equation}\label{equalphabeta}\begin{array}{cccc}
        (5.50,17.59)&(-24.20,-7.02)&(6.66,
        23.75)&(-30.30,-7.91)\\ \\ (9.39,
        18.75)&(-22.28,-10.91)&(9.56,
        23.02)&(-26.73,-10.76)\\ \\ 
        (11.90,19.50)&(-21.97,-13.25)&
        (11.61,22.71)&(-25.29,-12.71)
      \end{array}\end{equation}

    \begin{rems}\label{remscompare}
      {\rm (1) The definition of the quantities
        $\alpha_k$, $\beta_k$, $\gamma_k$ and $\delta_k$ implies that for
        $q\in (\bar{q}_7,\bar{q}_1]$, 
      $qx_{\natural}\in [0,23.76]$ and
      $\lim_{q\rightarrow (\bar{q}_7)^+}\in [0,23.76]$. Indeed,
      $23.76>\max_{1\leq k\leq 7}\beta_k$.

      (2) The table
        (\ref{equtable}) contains the first 8 double zeros of $\theta$. The
      maximal of their moduli is $4.08\ldots <23.76$.}
    \end{rems}   

  \item
We impose also the condition $qx_{\natural ,k}\in (x_{4k+3},x_{4k+5})~(*)$. 

This 
  is possible to obtain. Indeed, one has $qx_{4k+4}<x_{4k+3}<x_{4k+5}<qx_{4k+6}$,
  see (\ref{equineq}), and

  $$-q^{-4k-8}<x_{4k+8}<x_{4k+6}<-q^{-4k-6}<-q^{-4k-4}<x_{4k+4}<x_{4k+2}<-q^{-4k-2}~,$$
  so $qx_{4k+4}<-q^{-4k-3}<-q^{-4k-5}<qx_{4k+6}$. Part~(2) of Lemma~\ref{lm2zeros}
  implies 
  \begin{equation}\label{equ2ineq}
    |\theta (q,-q^{-4k-3})|<1~~~\, {\rm and}~~~\, |\theta (q,-q^{-4k-5})|<1~.
  \end{equation}
  We show
  that there exists $x_*\in (x_{4k+3},x_{4k+5})$ such that $x_*\in (-q^{-4k-3},-q^{-4k-5})$ and $\theta (q,x_*)>1$. Then we set $qx_{\natural ,k}:=x_*$.

  For $q\in [\bar{q}_{k+1},\bar{q}_k]$, the zeros $x_{4k+2}$ and $x_{4k+4}$
  are real and distinct. One can choose $x_{\dagger}\in (x_{4k+4},x_{4k+2})$ (hence $x_{\dagger}<0$ and 
  $\theta (q,x_{\dagger})<0$) such that
  $x_{\dagger}/q\in (x_{4k+3},x_{4k+5})$. Indeed, according to~(\ref{equFE}),

  \begin{equation}\label{equdagger}
    \theta (q,x_{\dagger}/q)=1+x_{\dagger}\theta (q,x_{\dagger})>1~.
    \end{equation}
  Thus $\theta (q,x_{\dagger}/q)>0$ and $x_{\dagger}/q\in (x_{4k+3},x_{4k+5})$.
  
  We prove that one has $-q^{-4k-3}<x_{\dagger}/q<-q^{-4k-5}~(**)$.
  Hence if we set
  $qx_{\natural ,k}:=x_{\dagger}/q$, i.~e. $x_*:=x_{\dagger}/q$, then we have $(*)$ and 

  \begin{equation}\label{equnatural}
    \theta (q,x_{\natural ,k})=1+qx_{\natural ,k}\theta (q,qx_{\natural ,k})>1~.
    \end{equation}
  Condition $(**)$ is fulfilled for $q$ close to $0$, if we choose
  $x_{\dagger}\sim -q^{-4k-3}$. Indeed, in this case
  $x_{4k+4}\sim -q^{-4k-4}$ and $x_{4k+2}\sim -q^{-4k-2}$, so
  $x_{\dagger}\in (x_{4k+4},x_{4k+2})$. Then one defines $x_{\dagger}$
  as continuously depending on $q\in [\bar{q}_{k+1},0)$. As
    conditions (\ref{equ2ineq}) are fulfilled and by (\ref{equdagger})
    $\theta (q,x_{\dagger}/q)>1$, 
    this means that by continuity on $q$ condition $(**)$ holds true for
    $q\in [\bar{q}_{k+1},0)$.
    
    \item 
      As we saw in the previous point, the quantity $x_{\natural, k}$ can be chosen
      to depend continuously on~$q$. Condition (\ref{equnatural}) and
      Lemma~\ref{lm-5} allow to conclude that $|\theta (q,z)|\neq 0$ for
      $|z|=|x_{\natural}|$, see~(\ref{equestimation}).
    
    \item 
      For $k=1$, $\ldots$, $5$ and for $q=\bar{q}_{k+1}$, one has
      $x_{\natural ,k}\in (-(\bar{q}_{k+1})^{-4k-6},-(\bar{q}_{k+1})^{-4k-4})$ and
      $x_{\natural ,k+1}\in (-(\bar{q}_{k+1})^{-4k-10},-(\bar{q}_{k+1})^{-4k-8})$. Hence
      $x_{\natural ,k+1}(\bar{q}_{k+1})>x_{\natural ,k}(\bar{q}_{k+1})$.
    \end{enumerate}

\subsection{The next spectral values\protect\label{subsecnextspectral}}

Suppose that $q\in (-1,\bar{q}_7]$.
  We remind that for $q\in (-1,0)$, the double zeros of $\theta$ belong to
  the interval $(-13.29,23.65)$, see \cite[Theorem~6]{KoSe}. We choose
  $qx_{\natural}$ (depending continuously on $q$) such that
  $qx_{\natural}\geq 24$ and $\theta (q,qx_{\natural})>0$.

  \begin{rem}\label{remxnatural}
    {\rm The function $qx_{\natural}$ has jump points at which it is continuous
      from the left; the left limit at a jump point is larger than the right
      limit. This is the case for $q\in (\bar{q}_7,\bar{q}_1]$, see point 4.
    in the previous subsection. The condition
    $qx_{\natural}\geq 24$ for $q\in (-1,\bar{q}_7]$ implies that at
  $\bar{q}_7$ this function also has such a jump point, see part (1) of
  Remarks~\ref{remscompare}. We choose the number $24$ for lower bound of
$qx_{\natural}$, because it is larger than $23.65$ and $23.76$.}
  \end{rem}
  Hence
  $x_{\natural}<0$ and $\theta (q,x_{\natural})>1$, see (\ref{equFE}).
  We denote by
  $I$ the interval of the form $[x_{4\nu -1},x_{4\nu +1}]$ with smallest
  possible value of $\nu$ which does not intersect the interval
  $[-13.29,24]$; 
  one has $\theta (q,.)\geq 0$ on~$I$, see
  Remark~\ref{remsigns}.
  Then we set $qx_{\natural}:=x_{4\nu -1}$.

  As $q$ decreases, both ends of $I$ tend to $1$ and for some value of $q$,
  one has $x_{4\nu -1}=24$. We denote by $q_{\diamond}$ the first, i.~e. the
  largest possible such
  value. For $q=q_{\diamond}$, one has to choose instead of
  $[x_{4\nu -1},x_{4\nu +1}]$ the interval $[x_{4\nu +3},x_{4\nu +5}]$ to play the
  role of the interval~$I$. Thus at the jump point $q_{\diamond}$,

  \begin{equation}\label{equsemicont}
    \lim_{q\rightarrow q_{\diamond}^-}qx_{\natural}>
    \lim_{q\rightarrow q_{\diamond}^+}qx_{\natural}~.
  \end{equation}
  We define $qx_{\natural}$ as semi-continuous from the left, i.~e. we set 
  $qx_{\natural}|_{q=q_{\diamond}}:=\lim_{q\rightarrow q_{\diamond}^-}qx_{\natural}$.

  The string of inequalities (\ref{equineq}) and Lemma~\ref{lm2zeros}
  imply that for $q=q_{\diamond}$,

  $$x_{4\nu +3}<qx_{4\nu +6}<-q^{-(4\nu +7)}~~~\, {\rm and}~~~\,
  -q^{-(4\nu -3)}<qx_{4\nu }<x_{4\nu -1}=24~,~~~\, {\rm so}$$

  $$\begin{array}{l}x_{4\nu +3}<-q_{\diamond}^{-(4\nu -3)}\cdot q_{\diamond}^{-10}<
    24\cdot q_{\diamond}^{-10}\\ \\
    <24\cdot (\bar{q}_7)^{-10}<24\cdot (0.91)^{-10}=61.630\ldots <61.64~.
    \end{array}$$
Theorem~\ref{tmlocation} is proved.

\section{Proof of Theorem~\protect\ref{tmimprove}\protect\label{sectmimprove}}

\subsection{Proof of part (1)}


  We prove a proposition from which part (1) of the theorem follows:

\begin{prop}\label{propRe2.4}
  For $q\in (0,1)$ and {\rm Re}$x=2.4$, the function $\theta$ has no zeros.
\end{prop}
The zeros of $\theta$ depend continuously on $q$ and for $q\leq \tilde{q}_1$,
all zeros are real negative. The absence of zeros with Re$x=2.4$ implies
their absence for Re$x\geq 2.4$. Indeed, no zeros are born at infinity.

\begin{proof}[Proof of Proposition~\ref{propRe2.4}] 
  For $q\in (0,q_{\heartsuit}]$, the function $\theta$
      has no zeros with non-negative real part,
      see Proposition~\ref{propheart}.
     For $x>0$ fixed and for $q\in [q_{\heartsuit},1)$,
        the minimal possible value of
        $\theta$ is obtained for $q=q_{\heartsuit}$. Indeed, all terms of the
        series of
        $\theta$ are positive. A direct computation yields
        $\theta (q_{\heartsuit},2.4)=9.1023\ldots$.

        For $q\in (0,1)$ and Re$x\geq 2.4$, one has
        $|G(q,x)|\leq \sum_{j=1}^{\infty}(1/2.4^j)=1/1.4=0.71428\ldots =:g_0$.
        Set $x:=a+bi$. For $b=0$, and for $a>0$ and $q\in [q_{\heartsuit},1)$
          fixed,
        all terms of the series of $-G$
        are positive. As $b$ increases, their moduli become smaller and
        their arguments are no longer all of them the same, so $|G|$ decreases.
        Thus

        $$|G|\leq g_0~~~\, {\rm for}~~~\, 
        (a,b,q)\in \mathcal{W}:=
        [2.4,\infty )\times [0,\infty )\times [q_{\heartsuit},1)~.$$
        We give a lower bound for $|\Theta^*|$ for $(a,b,q)\in \mathcal{W}$. 
        First of all we notice that
        $|1+q^k/x|=|1+q^ka/(a^2+b^2)-bi/(a^2+b^2)|>1$, so

        $$|\Theta^*|>(\prod_{k=1}^{\infty}(1-q^k))\cdot
        (\prod_{k=1}^{\infty}|1+q^kx|)~.$$
        For $b\geq 0$ and $(a,q)$ fixed, the modulus $|1+q^kx|=|1+q^ka+q^kbi|$
        is minimal for $b=0$. Therefore

        $$|\Theta^*|>\prod_{k=1}^{\infty}((1-q^k)(1+q^ka))\geq
        \prod_{k=1}^{\infty}((1-q^k)(1+q^k\cdot 2.4))=:T~.$$
        We prove that $\ln T>0$ hence $|\Theta^*|>T>1>g_0$ and
        $\theta =\Theta^*-G\neq 0$ for Re$x=2.4$. 
              
        We set $q:=e^{-t}$, $t\in [0,-\ln q_{\heartsuit}]$, $-\ln q_{\heartsuit}=0.3236567932\ldots $,
        
          $$U:=r\sum_{k=1}^{\infty}\ln (1-e^{-kr})~~~\, {\rm and}~~~\, 
          V:=r\sum_{k=1}^{\infty}\ln (1+2.4\cdot e^{-kr})~.$$
          Our first observation is that the functions $f:=\ln (1-e^{-t})$ and
          $h:=\ln (1+2.4\cdot e^{-t})$ are monotone, increasing from $-\infty$ to $0$
          and decreasing
          from $\ln 3.4$ to $0$ respectively. The quantities $U$ and $V$
          are their Riemann's sums. Next, we find that

          $$\int_0^{\infty}f(t)dt=-1.644933967\ldots =:I_f~~~\, {\rm and}~~~\,
          \int_0^{\infty}h(t)dt=1.648348390\ldots =:I_h>|I_f|~.$$ 
          The monotonicity of the functions $f$ and $h$ implies that the
          quantities $|U|$ and $V$ are decreasing as $r$
          increases in $[0,-\ln q_{\heartsuit}]$. A numerical check
          shows that for $r\in [0.005,-\ln q_{\heartsuit}]$, one has $V>U+0.01$
          (it suffices to
          check this for the truncations of $U$ and $V$ up to the 2000th term).
          For $r=0.005$, one obtains $V=1.645\ldots >|I_f|$. Hence it is true
          that for $r\in (0,0.005]$, $I_h>V(r)>|I_f|>|U(r)|$. Thus
          $\sum_{k=1}^{\infty}|\ln (1-e^{-kr})|<\sum_{k=1}^{\infty}\ln (1+2.4\cdot e^{-kr})$
          and $\ln |\Theta^*|>0$. Therefore $|\Theta^*|>1>g_0\geq |G|$ and
          $\theta \neq 0$. 
          
  \end{proof}

\subsection{Proof of part (2)}

We deduce the proof of part (2) from several lemmas whose proofs are given
in the next subsection. In the formulations of parts (2) of Lemmas~\ref{lmA} and \ref{lmB} we use the interval $[\pi /3,\pi /2]$ instead of the minimal necessary $[\pi /3, 5\pi /12]$ for technical reasons. In the proofs of the lemmas we sometimes give estimations not on the minimal necessary intervals for~$q$, but on larger ones.  

\begin{lm}\label{lmA}
  (1) For $x=1.7e^{\phi i}$, $\phi \in [5\pi /12,\pi /2]$, and
  $q\in [q_{\heartsuit},0.88]$,
  one has $\theta (q,x)\neq 0$.
  \vspace{1mm}
  
(2) For $x=1.5e^{\phi i}$,
  $\phi \in [\pi /3,\pi /2]$, and $q\in [q_{\heartsuit},0.88]$, one has $\theta (q,x)\neq 0$.
  \vspace{1mm}
  
  (3) For $x=1.4e^{\phi i}$, $\phi \in [\pi /4,\pi /3]$ and
  $q\in [q_{\heartsuit},0.88]$,
  one has $\theta (q,x)\neq 0$.

  \end{lm}

\begin{lm}\label{lmB}
(1) For $x=1.7e^{\phi i}$, $\phi \in [5\pi /12,\pi /2]$, and $q\in [0.88,1)$,
    one has $|\Theta^*(q,x)|<0.021$.

    (2) For $x=1.5e^{\phi i}$, $\phi \in [\pi /3,\pi /2]$, and $q\in [0.88,1)$,  one
    has $|\Theta^*(q,x)|<0.1537$.
    
 (3) For $x=1.4e^{\phi i}$, $\phi \in [\pi /4,\pi /3]$, and $q\in [0.88,1)$,
    one has $|\Theta^*(q,x)|<0.5087$.
    
  \end{lm}

\begin{lm}\label{lmC}
  (1) For $x=1.7e^{\phi i}$, $\phi \in [5\pi /12,\pi /2]$, and $q\in [0.88,1)$,
    one has
  $|G|>0.34689$.

  (2) For $x=1.5e^{\phi i}$, $\phi \in [\pi /3,5\pi /12]$, and $q\in [0.88,1)$, one has
  $|G|>0.4402$.

  (3) For $x=1.4e^{\phi i}$, $\phi \in [\pi /4,\pi /3]$, and $q\in [0.88,1)$,
    one has $|G|>0.5960$.
  \end{lm}

\begin{lm}\label{lmD}
(1) The function $\theta$ has no zero on either of the two segments
  $[1.5e^{\pm 5\pi i/12},1.7e^{\pm 5\pi i/12}]$, for any $q\in (0,1)$.

  (2) The function $\theta$ has no zero on either of the two segments
  $[0,1.5e^{\pm \pi i/3}]$, for any $q\in (0,1)$.

  (3) The function $\theta$ has no zero on either of the two segments
  $[0,1.4e^{\pm \pi i/4}]$, for any $q\in (0,1)$.

  (4) For $q\in (0,1)$, the function $\theta$ has no zeros on the segment
  $[-1.7i,1.7i]$.
  \end{lm}

Lemma~\ref{lmD} shows that there are no zeros of $\theta$ on any of the
segments which are parts of the boundary $\partial \tilde{R}$
of the domain $\tilde{R}$. 
From Lemma~\ref{lmA} one deduces the absence of zeros on the arcs
of $\partial \tilde{R}$ for $q\in [q_{\heartsuit},0.88]$; 
there are no zeros of $\theta$ in the right half-plane
for $0<q<q_{\heartsuit}$ (Proposition~\ref{propheart}).
Lemmas~\ref{lmB} and~\ref{lmC} imply the absence of zeros on these arcs for
$q\in [0.88,1)$.

For $q\in (0,\tilde{q}_1]$, all zeros of $\theta$ are real negative. Zeros
depend continuously on the parameter~$q$. The absence of zeros on
$\partial \tilde{R}$ implies that no zero enters the domain $\tilde{R}$,
i.~e. for all $q\in (0,1)$, $\tilde{R}$ is free of zeros of~$\theta$.

\subsection{Proofs of the lemmas}

\begin{proof}[Proof of Lemma~\ref{lmA}]

Part (1). Set $u:=\cos (\phi )$ (so $\sin \phi =\sqrt{1-u^2}$),
  $T:=\sum_{j=0}^9q^{j(j+1)/2}x^j$ and $T_I:=$Im$T$. 
  Then $T_I$ is of the form
  $\sqrt{1-u^2}T_I^0$, where $T_I^0:=\sum_{j=0}^8s_ju^j$, $s_j\in \mathbb{R}[q]$.

  The coefficients $s_0$, $s_3$, $s_4$, $s_7$ and $s_8$ are positive-valued on
  $[q_{\heartsuit},0.88]$.
  Moreover, one has $s_0>0.47
  $. It is true that
  
  $$\sigma_0:=\sum_{j=10}^{\infty}0.88^{j(j+1)/2}1.7^j=0.2918\ldots ~~~\,
  {\rm and}~~~\, \min_{\phi \in [5\pi /12, \pi /2]}\sin (\phi )=\sin (5\pi /12)=0.965925\ldots ~.$$
Suppose first that $q\in [q_{\heartsuit},0.82]$. All coefficients $s_j$ are positive-valued, so 

$$|\theta (q,x)|\geq |T(q,x)|-\sum_{j=10}^{\infty}q^{j(j+1)/2}|x^j|\geq
|T_I|-\sigma_0\geq \sin (5\pi /12)\cdot s_0-\sigma_0>0$$
hence $\theta (q,x)\neq 0$. Now suppose that $q\in [0.82,0.88]$.
We introduce the following notation:

  $$\begin{array}{ccl}v(\alpha )&:=&\alpha \cdot |s_1|+\alpha^2\cdot |s_2|
    +\alpha^5\cdot |s_5|+\alpha^6\cdot |s_6|~,\\ \\
    w(\beta )&:=&\beta^3\cdot s_3+\beta^4\cdot s_4+
    \beta^7\cdot s_7+\beta^8s_8~.\end{array}$$
For $u\in [0,0.14]$, the sum $s^{(1)}:=s_0-v(0.14)$ 
  is $>0.307$. On the other hand $T_I^0\geq s^{(1)}$. The following inequalities
  hold true:

  \begin{equation}\label{equthetas}\begin{array}{ccccccc}|\theta (q,x)|&\geq&
    |T(q,x)|-\sum_{j=10}^{\infty}q^{j(j+1)/2}|x^j|&\geq&
    |T_I|-\sigma_0\\ \\ &\geq &\sin (5\pi /12)\cdot  |T_I^0|-0.2919&\geq&
    \sin (5\pi /12)\cdot  s^{(1)}-0.2919&>&0~,
    \end{array}\end{equation}
  therefore $\theta (q,x)\neq 0$. We define a finite sequence of numbers
  $a_j$, $j=1$, $\ldots$, $14$:
  
  $$\begin{array}{lllllll}0.14~,& 0.15~,& 0.16~,& 0.17~,& 0.18~,& 0.187~,& 0.194~,\\ \\
    0.201~,& 0.208~,& 0.22~,& 0.23~,& 0.24~,& 0.25~,&\cos (5\pi /12)=0.258819\ldots ~.
     \end{array}
$$
  We check numerically that $s^{(j+1)}:=s_0-v(a_{j+1})+w(a_j)>0.34$, $j=1$, $\ldots$, $13$.
  The sum $s^{(j+1)}$ is a minoration of $T_I^0$ for $u\in [a_j,a_{j+1}]$.
  By analogy with (\ref{equthetas}) we conclude that $\theta (q,x)\neq 0$. 
\vspace{1mm}
  
Part (2). We set $x:=1.5(\cos (\phi )+i\sin (\phi ))$,
  $\phi \in [\pi /3,\pi /2]$, and $T:=\sum_{j=0}^9q^{j(j+1)/2}x^j=
\tilde{T}_R+i\tilde{T}_I$, where $\tilde{T}_R:=$Re$T$ and $\tilde{T}_I:=$Im$T$.
Suppose first that $q\in [q_{\heartsuit},0.84]$. 
%
  Set $u:=\cos \phi$, so
  $\sin \phi =\sqrt{1-u^2}$. In the variables $(q,u)$ we write $T_R$, $T_I$
  instead of $\tilde{T}_R$, $\tilde{T}_I$. Then $T_I$ is of the form
  $\sqrt{1-u^2}T_I^0$, $T_I^0\in \mathbb{R}[q,u]$. 
  
  Set
  $S(q,u):=\partial T_I^0/\partial u$. Hence $S=\sum_{j=0}^7S_ju^j$,
  $S_j\in \mathbb{R}[q]$. For $q\in [q_{\heartsuit},0.84]$, all coefficients $S_j$
  are positive-valued, so $T_I^0$ is minimal for $u=0$.
  For $q\in [q_{\heartsuit},0.88]$, the polynomial $T_I^0|_{u=0}$
  takes only values $>0.5$. The sum $T_{\bullet}:=\sum_{j=10}^{\infty}0.88^{j(j+1)/2}1.5^j$
  is $<0.0781$ and one has 
  $\sqrt{3}/2\leq \sqrt{1-u^2}\leq 1$. Thus

  \begin{equation}\label{equTT}
    |\theta (q,x)|>|T_I|-T_{\bullet}\geq T_I^0|_{u=0}\cdot \sqrt{3}/2-0.0781>0.5\cdot \sqrt{3}/2-0.0781>
    0.4330-0.0781>0~.
    \end{equation}
  Suppose now that $q\in [0.84,0.88]$. The coefficients $S_0$, $S_1$, $S_4$
  and $S_5$ are decreasing while $S_2$, $S_3$, $S_6$ and $S_7$ are increasing.
  We define the finite sequence of numbers $b_j$, $1\leq j\leq 8$:
  
  $$0.84~,~~~\, 0.85~,~~~\, 0.86~,~~~\, 0.868~,~~~\, 0.874~,~~~\, 0.878~,~~~\, 0.8795~,~~~\, 0.88~.$$
  For $q\in [b_j,b_{j+1}]$, $1\leq j\leq 7$, the function $S$ is bounded from below by the polynomial 
  
  $$S_0(b_{j+1})+S_1(b_{j+1})u+S_2(b_j)u^2+S_3(b_j)u^3+S_4(b_{j+1})u^4+S_5(b_{j+1})u^5+S_6(b_j)u^6+S_7(b_j)u^7~.$$
  All these $7$ polynomials are positive-valued for $u\in [0,1/2]$, so $S>0$
  and by (\ref{equTT}), $|\theta (q,x)|>0$, so $\theta (q,x)\neq 0$. 
  \vspace{1mm}
  
Part (3). We consider the quantities $T:=\sum_{j=0}^9q^{j(j+1)/2}x^j$ and
  $T_I:=$Im$T$.
  Setting $u:=\cos \phi$ (hence $\sin \phi =\sqrt{1-u^2}$) one can represent
  $T_I$ in the form $\sqrt{1-u^2}T_I^0$, where $T_I^0\in \mathbb{R}[u,q]$
  is a degree $8$ polynomial in~$u$. Set $S:=\partial T_I^0/\partial u$.
  Then $S=\sum_{j=0}^7s_ju^j$, $s_j\in \mathbb{R}[q]$. One finds numerically
  that for $q\in [q_{\heartsuit},0.84]$, all coefficients $s_j$ are positive-valued and
  $s_0>0.5$. Hence $S>0$ and $T_I^0$ is minimal for $u=1/2$, with $T_I^0|_{u=1/2}>1.1$;
  the latter lower bound is valid for $q\in [q_{\heartsuit},0.88]$.
  The sum $T_{\dagger}:=\sum_{j=10}^{\infty}0.88^{j(j+1)/2}1.4^j$ is $<0.038$, so 
  
  $$|\theta (q,x)|\geq |T_I|-T_{\dagger}> \min_{u\in [1/2,\sqrt{2}/2]}\sqrt{1-u^2}\cdot
  T_I^0|_{u=1/2}-0.038>(\sqrt{2}/2)\cdot 1.1-0.038>0~.$$
  Suppose that $q\in [0.84,0.88]$. Only the coefficients $s_4$ and $s_5$ can
  take negative values, the other coefficients are positive-valued.
  Hence for $u\in [u_1,u_2]$, the following function is a minoration of~$S$:
  
  $$S^{\bullet \bullet}:=\sum_{j\in J}s_j\cdot (u_1)^j-|s_4|\cdot (u_2)^4
  -|s_5|\cdot (u_2)^5~,~~~\, J:=\{ 0,1,2,3,6,7\}~.$$
  One finds that for the following values of the couple $(u_1,u_2)$,
  one has $S^{\bullet \bullet}>0.187$:
  
  $$(0.5,0.55)~,~~~\, (0.55,0.6)~,~~~\, (0.6,0.65)~,~~~\,
  (0.65,\sqrt{2}/2)~~~\, ({\rm with}~\,  \sqrt{2}/2=0.7071\ldots )~.$$
  Hence as above $|\theta (q,x)|>0$, i.~e. $\theta (q,x)\neq 0$.

  \end{proof}

\begin{proof}[Proof of Lemma~\ref{lmB}]
  In the proof of the lemma we use the following technical result:

  \begin{lm}\label{lmtech}
  For $x=Ae^{\phi i}$, $A>0$, $\phi \in [\alpha ,\beta ]\subset [0,\pi /2]$,
  and for $q\in (0,1)$, $k\in \mathbb{N}$ fixed,
  the quantity $r_k:=|(1+q^kx)(1+q^{k-1}/x)|$ is maximal for
$\phi =\alpha$. 
  \end{lm}

\begin{proof}
  Set $u:=\cos \phi$ (hence $\sin \phi =\sqrt{1-u^2}$) and $Q:=q^{k-1}$.
  One finds that
Re$r_k=(A^2Qqu+AQ^2q+Qu+A)/A$ and Im$r_k=Q(A^2q-1)\sqrt{1-u^2}/A$. Thus 
  
$$r_k^2=4Q^2qu^2+Ku+L~,~~~\, K:=2AQ^3q^2+2AQq+2Q^3q/A+2Q/A~,~~~\,
L:=A^2Q^2q^2+Q^2/A^2+q^2Q^4+1$$
and hence $\partial r_k^2/\partial u=8Q^2qu+K>0$. This means that $r_k$
  is maximal when $u$ is maximal, i.~e. for $u=\cos \alpha$. 
  \end{proof}

We prove part (2) of Lemma~\ref{lmB} first. We apply Lemma~\ref{lmtech}
to the case
$[\alpha ,\beta ]=[\pi /3,\pi /2]$,
$A=1.5$. For fixed $q$, the quantity $|\Theta^*|$ is maximal for
$x=1.5e^{\pi i/3}$. 
We set

$$\xi :=1.5e^{\pi i/3}~~~\, {\rm and}~~~\, f(q):=(1-q)(1+q\xi )(1+q/\xi )~.$$
Thus
$\Theta^*(q,\xi )=(1+1/\xi )\prod_{k=1}^{\infty}f(q^k)$. The function
$|f|\in C_{[0,1]}$ has a single critical point in $(0,1)$
at $q^{\circ}:=0.1841745227\ldots$ which is a maximum, with 
$|f(q^{\circ})|=f^{\circ}:=1.012099679\ldots$; $|f|$ is increasing on
$[0,q^{\circ}]$ and decreasing on $[q^{\circ},1]$. It takes the value $1$
only at $0$ and at $q^{\bullet}:=0.3100969623\ldots$.

We give an upper bound for the quantity $|\Theta^*(q,\xi )|$ for
$q\in [0.88,1)$. 

\begin{lm}
 For $q\in [0.88,1)$, one has $|\Theta^*(q,\xi )|<0.1537$.
\end{lm}

\begin{proof}
We set $q:=e^{-t}$ and 
we consider the function 

$$l(e^{-t},\xi ):=\ln |\Theta^*(e^{-t},\xi )|=
\ln |1+1/\xi |+l_0(e^{-t},\xi )~,~~~\,
l_0(e^{-t},\xi ):=\sum _{k=1}^{\infty}\ln |f(e^{-kt})|~.$$
The function $\ln |f(e^{-t})|$

\begin{itemize}
  \item is increasing for $t\in (0,-\ln (q^{\circ})]$
 (with $-\ln (q^{\circ}):=t^{\sharp}=1.691871478\ldots$);
 
 \item is decreasing for
   $t\in [-\ln (q^{\circ}),\infty )$;

   \item equals $0$ for
     $t=t^{\bullet}:=-\ln (q^{\bullet})=1.170870249\ldots$;

     \item tends to $-\infty$ as
       $t\rightarrow 0^-$ and to $0^+$ as $t\rightarrow \infty$;

       \item the integral
         $I:=\int_{t=t^{\bullet}}^{\infty}\ln |f(e^{-t})|dt$ converges.
         \end{itemize}

  When $q\in [0.88,1)$, then $-\ln q\in (0 ,t^{\flat}]$,
  $t^{\flat}:=-\ln (0.88)=0.1278333715\ldots$. We give first an upper bound of
  $\ln |\Theta^*(e^{-t^{\flat}},\xi )|$. In the sum
  $l_0(e^{-t^{\flat}},\xi )$, negative are the terms with $k=1$, $\ldots$, $9$,
  and their
  sum is $-2.4132\ldots$. The sum of the positive terms multiplied by
  $t^{\flat}$ 
  is

  $$\leq I+\ln |f(q^{\circ})|\cdot t^{\flat}=0.01972190427\ldots
  +0.01202706294\ldots \cdot t^{\flat}=
  0.02125\ldots <0.02126~.$$
  Indeed, the sum of the positive terms with values of $k$ such that
  $kt^{\flat}\leq t^{\sharp}$ (resp. $kt^{\flat}>t^{\sharp}$)
  is equal to the sum of the surfaces of the
  rectangles with bases $[kt^{\flat},(k+1)t^{\flat}]$ (resp.
  $[(k-1)t^{\flat},kt^{\flat}]$) and height
  $\ln |f(e^{-kt^{\flat}})|$. The term $\ln |f(q^{\circ})|\cdot t^{\flat}$
  is added to compensate
  the two overlapping rectangles with bases containing $t^{\sharp}$.

  Thus the sum $l_0(e^{-t^{\flat}},\xi )$ is $\leq -2.4132+0.02126/t^{\flat}=
  -2.24688975\ldots 
  <-2.2468$, so

  $$\ln |\Theta^*(e^{-t^{\flat}},\xi )|\leq \ln |1+1/\xi |-2.2468=
  0.3736072006\ldots -2.2468<-1.8731~.$$
  One finds that $|\Theta^*(q^{\circ},\xi )|\leq e^{-1.8731}=0.1536\ldots <0.1537$.
  However we need an estimation for $q\in [0.88,1)$, with $0.88>q^{\circ}$.

  If $q\in (q^{\circ},1)$ (hence $t\in (0,t^{\flat})$), then the above estimation
  of the sum of the positive terms in $l_0(e^{-t},\xi )$ remains exactly the
  same while the negative terms become smaller and eventually new negative
  terms appear. This is why the inequality $|\Theta^*(q,\xi )|<0.1537$
  remains true.
  \end{proof}

  We prove part (3) next. Making use of  Lemma~\ref{lmtech}
with 
$[\alpha ,\beta ]=[\pi /4,\pi /3]$ and 
$A=1.4$, we observe that for fixed $q$, the quantity $|\Theta^*|$ is maximal for
$x=1.4e^{\pi i/4}$, so we set $\eta :=1.4e^{\pi i/4}$ and
  $f_1(q):=(1-q)(1+q\eta )(1+q/\eta )$. Thus
$\Theta^*(q,\eta )=(1+1/\eta )\prod_{k=1}^{\infty}f_1(q^k)$. The function
$|f_1|\in C_{[0,1]}$ has a single critical point in $(0,1)$
at $q^{\spadesuit}:=0.2807210416\ldots$ which is a maximum and 
$|f_1(q^{\spadesuit})|=f_1^{\spadesuit}:=1.082269345\ldots$; $|f_1|$ is increasing on
$[0,q^{\spadesuit}]$ and decreasing on $[q^{\spadesuit},1]$. It takes the value $1$ only 
  at $0$ and at $q^{\bullet \bullet}:=0.5079981947\ldots$.

\begin{lm}
 For $q\in [0.88,1)$, one has $|\Theta^*(q,\eta )|<0.5087$.
\end{lm}

\begin{proof}
We set $q:=e^{-t}$ and 
we consider the function 

$$l(e^{-t},\eta ):=\ln |\Theta^*(e^{-t},\eta )|=
\ln |1+1/\eta |+l_1(e^{-t},\eta )~,~~~\,
l_1(e^{-t},\eta ):=\sum _{k=1}^{\infty}\ln |f_1(e^{-kt})|~.$$
The function $\ln |f_1(e^{-t})|$
\begin{itemize}
\item
  is increasing for $t\in (0,-\ln (q^{\spadesuit})]$
 (with $-\ln (q^{\spadesuit}):=t^{\spadesuit}=1.270393837\ldots$);
\item
  is decreasing for
   $t\in [-\ln (q^{\spadesuit}),\infty )$;
  \item
    equals $0$ for
    $t=t^{\bullet \bullet}:=-\ln (q^{\bullet \bullet})=0.6772773850\ldots$;
  \item
    tends to
   $-\infty$ as
    $t\rightarrow 0^-$ and to $0^+$ as $t\rightarrow \infty$;
  \item
    the integral
  $I_1:=\int_{t=t^{\bullet \bullet}}^{\infty}\ln |f_1(e^{-t})|dt$ converges.
\end{itemize}

  When $q\in [0.88,1)$, then $-\ln q\in (0 ,t^{\flat}]$,
  $t^{\flat}:=-\ln (0.88)=0.1278333715\ldots$. We give an upper bound of
  $\ln |\Theta^*(e^{-t^{\flat}},\eta )|$. In the sum
  $l_1(e^{-t^{\flat}},\eta )$, negative are the terms with $k=1$, $\ldots$, $5$,
  and their
  sum is $-3.567495877\ldots$. The sum of the positive terms multiplied by
  $t^{\flat}$ 
  is

  $$\leq I_1+\ln |f_1(q^{\spadesuit})|\cdot t^{\flat}=0.1464923181\ldots
  +0.1581201638\ldots \cdot t^{\flat}=
  0.1667053517\ldots <0.166706~.$$

  Thus the sum $l_1(e^{-t^{\flat}},\eta )$ is $\leq -3.567495+0.166706/t^{\flat}=
  -2.263406732\ldots 
  <-2.2634$, so

  $$\ln |\Theta^*(e^{-t^{\flat}},\eta )|\leq \ln |1+1/\eta |-2.2634=
  1.587563109\ldots -2.2634<-0.6758~.$$
  One finds that $|\Theta^*(q^{\spadesuit},\eta )|\leq e^{-0.6758}=
  0.5087492584\ldots <0.5087$.

  If $q\in (q^{\spadesuit},1)$ (hence $t\in (0,t^{\flat})$), then the above
  estimation
  of the sum of the positive terms in $l_1(e^{-t},\eta )$ remains exactly the
  same while the negative terms become smaller and eventually new negative
  terms appear. This is why the inequality $|\Theta^*(q,\eta )|<0.5087$
  remains true.
  \vspace{1mm}

  Finally we prove part (1). We use Lemma~\ref{lmtech}
with 
$[\alpha ,\beta ]=[5\pi /12,\pi /2]$,
$A=1.7$. We set $x:=1.7e^{\phi i}$,
  $\phi \in [5\pi /12,\pi /2]$. As $|\Theta^*|$ is maximal for
  $\phi =5\pi  /12$, we set $\zeta :=1.7e^{5\pi i/12}$ and we consider the
  function $f_2(q):=(1-q)(1+q\zeta )(1+q/\zeta)$. One has
  $\Theta^*(q,\zeta )=(1+1/\zeta )\prod_{k=1}^{\infty}f_2(q^k)$.

  The function $|f_2(q)|$ is decreasing from $1$ to $0$.
  As in the proofs of parts (2) and (3)
  we define the functions

  $$l(e^{-t},\zeta ):=\ln |\Theta^*(e^{-t},\zeta )|=\ln |1+1/\zeta |+
  l_2(e^{-t},\zeta )~~~\, {\rm  and}~~~\,
  l_2(e^{-t},\zeta )=\sum_{k=1}^{\infty}\ln |f_2(e^{-kt})|~.$$
  The function $\ln |f_2(e^{-t})|$ is increasing from $-\infty$ to $0$.
  Suppose that $q\in [0.88,1)$. Then $t\in (0,t^{\flat}]$. The negative sum
  $l_2(e^{-t},\zeta )$ is maximal for $t=t^{\flat}$ when it equals
  $-4.153\ldots <-4.153$. Thus

  $$l(e^{-t},\zeta )<\ln |1+1/\zeta |-4.153=0.25054\ldots -4.153<-3.90245$$
  and $|\Theta^*|<e^{-3.90245}=0.020\ldots <0.021$.

  \end{proof}

  \end{proof}

\begin{proof}[Proof of Lemma~\ref{lmC}]

We need another lemma:

\begin{lm}\label{lmG1}
  Set $y:=r(u-\sqrt{1-u^2}i)$, $G_1:=y+qy^2+q^3y^3$ and $G_0:=y+qy^2+q^3y^3+q^6y^4$.
  Suppose that $(q, r)\in [0,1]^2$. Then 
 \vspace{1mm}
 
 (1) for $u\in [0,1]$, one has $\partial |G_1|^2/\partial r>0$;
 
 (2) for $u\in [1/4,1]$ (with $1/4<\cos (5\pi /12)$), one has $\partial |G_0|^2/\partial r>0$. 
\end{lm}

\begin{proof}
  As $\partial |G_1|^2/\partial r=2r(G_2+G_3)$, where
  $G_2=(5q^4r^3+3qr)u+8q^3r^2u^2\geq 0$ and 

$$G_3=3q^6r^4-2q^3r^2+1+2r^2(q^2-q^3)\geq 3q^6r^4-2q^3r^2+1>0$$
(the discriminant of the quadratic polynomial in $q^3r^2$ is negative),
one has $\partial |G_1|^2/\partial r>0$. Next, 
$\partial |G_0|^2/\partial r=2r(G_4+G_5)$, where
$G_4=2r^2(q^2-q^3)+5r^3u(q^4-q^6)+3r^4(q^6-q^7)\geq 0$ and 

$$G_5=20q^6r^3u^3+(12q^4r^2+8)q^3r^2u^2+(7q^8r^4-10q^5r^2+3)qru+4q^{12}r^6-3q^7r^4-2q^3r^2+1~.$$
We prove first that $g^*:=G_5|_{u=1/4}>0$. Indeed, 

$$16g^*=64q^{12}r^6+28q^9r^5-36q^7r^4-35q^6r^3-24q^3r^2+12qr+16~.$$
The inequality $16g^*>0$ results from 

$$\begin{array}{l}
  36q^{12}r^6-36q^7r^4+9qr\geq 36q^{12}r^6-36q^7r^4+9q^2r^2=9q^2r^2(2q^5r^2-1)^2\geq 0~,\\ \\
  28q^{12}r^6-28q^6r^3+7=7(2q^6r^3-1)^2\geq 0~,\\ \\
  7q^9r^5-7q^6r^3+7qr/4\geq (7/4)q^3r^2(4q^6r^2-4q^3r+1)=
  (7/4)q^3r^2(2q^3r-1)^2\geq 0~~~\, {\rm and}\\ \\
  \tilde{\varphi }(q,r):=21q^9r^5-24q^3r^2+5qr/4+9>0~~~\, (A).
  \end{array}$$
Of these inequalities only $(A)$ is not evident.
We need to prove it only for $q\geq 0.6$.
Indeed, the function $\tilde{h}:=-24x^2+5x/4+9$ takes only positive values
for $x\in [0,0.6]$ and $\tilde{\varphi}(q,r)\geq \tilde{h}(qr)$.
For $q\geq 0.6$ the proof of inequality $(A)$ is numerical.
Namely, we list the finite strictly increasing sequence denoted by
$v_j$, $1\leq j\leq 42$; we set $v_0:=0.6$:

$$\begin{array}{lllllllllll}
  0.7&0.76&0.78&0.79&0.8&0.81&0.815&0.82&0.825&0.828&0.831\\ \\
  0.834&0.837&0.84& 0.843&0.846&0.849&0.852&0.855&0.858&0.861&0.864\\ \\
  0.867&0.87&0.875&0.88&0.886&0.892& 0.898&0.904&0.911&0.918&0.925\\ \\
  0.933&0.941&0.95&0.96&0.97&0.98&0.99&1&1.01~.
  \end{array}$$
We check that each of the functions
$\tilde{\varphi}_j:=21v_j^9r^5-24v_{j+1}^3r^2+5v_jr/4+9$, $0\leq j\leq 41$,
takes only positive values on the interval $[0,1]$.
The function $\tilde{\varphi}_j$ is a minoration of $\tilde{\varphi}$
for $q\in [v_j,v_{j+1}]$.
Thus $g^*>0$. Now we show that for $u\in [1/4,1]$, $\partial G_5/\partial u>0$
from which $G_5>0$ follows. It is clear that $\partial^2G_5/\partial u^2>0$.
Hence $\partial G_5/\partial u$ is minimal for $u=1/4$ when it equals
$(qr/4)(28q^8r^4+24q^6r^3-25q^5r^2+16q^2r+12)>0$
(because $16q^5r^2\leq 16q^2r$ and $9q^5r^2\leq 9$).
\end{proof}

Part (1) of Lemma~\ref{lmC}. We  set $y:=1.7^{-1}(u-\sqrt{1-u^2}i)$
(hence $y=1/x$),  
  $u=\cos \phi \in [0,\cos (5\pi /12)]$, $\cos (5\pi /12)=0.258819\ldots$
  and as above $G_1:=y+qy^2+q^3y^3$. The quantity $|G_1|^2\in \mathbb{R}[q,u]$ is a
  degree 3 polynomial in $u$ with coefficients in $\mathbb{R}[q]$,
  positive-valued for $q\in [0.88,1]$. Thus
  $|G_1(q,0)|\leq |G_1(q,u)|\leq |G_1(q,\cos (5\pi /12))|$. The modulus
  $|G_1(q,0)|$
  is decreasing while $|G_1(q,\cos (5\pi /12))|$ is increasing in $q$. Hence

  $$0.51672<0.516720\ldots\leq |G_1(0.88,0)|\leq
  |G_1(q,u)|\leq |G_1(1,\cos (5\pi /12))|=0.66454\ldots <0.66455~.$$
  It is true that $G(q,y)=G_1(q,y)+G^{\diamond}$,
  $G^{\diamond}:=\sum_{j=1}^{\infty}q^{3j(3j-1)/2}y^{3j}G_1(q,q^{3j}y)$.
  From Lemma~\ref{lmG1} we deduce that $|G_1(q,q^{3j}y)|\leq |G_1(q,y)|$.
  This implies

  $$\begin{array}{l}
    |G^{\diamond}|\leq \sum_{j=1}^{\infty}1.7^{-3j}\cdot 0.66455=:g_{\diamond}=
    0.1698298885\ldots~~~\,
    {\rm and}\\ \\ 
  |G|\geq |G_1(q,y)|-g_{\diamond}\geq 0.51672-g_{\diamond}>0.34689~.\end{array}$$
\vspace{1mm}
  
Part (2) of Lemma~\ref{lmC}. We set $y:=(2/3)(u-\sqrt{1-u^2}i)$, $u=\cos \phi \in [0,1/2]$,
  so $y=1/x$, and as above
   $G_0:=y+qy^2+q^3y^3+q^6y^4$. Then $|G_0|^2\in \mathbb{R}[u,q]$ is a
  degree 3 polynomial in $u$ with coefficients in $\mathbb{R}[q]$. One has

  $$\partial (|G_0|^2)/\partial u=256q^6u^2/81+(512q^7/729+128q^3/81)u+
  (256q^9-1728q^6+576q^4+1296q)/2187~.$$
  This polynomial is positive-valued for $u\in [0,1/2]$, $q\in (0,1]$. Indeed,
  $\partial^2 (|G_0|^2)/\partial u^2>0$ and $\partial |G_0|^2/\partial u|_{u=0}$
  is positive-valued for $q\in (0,1]$, so $\partial |G_0|^2/\partial u>0$. Hence $|G_0|^2$ is minimal
  for $u=0$ and maximal for $u=1/2$.

  One finds numerically that the minimal value of $|G_0|^2|_{u=0}$ is
  $0.4451\ldots >0.4451=:g_0$ (attained for $q=1$) and that for $q\in (0,1]$,
  $|G_0|^2|_{u=1/2}<0.85$. On the other hand

  $$G(q,y)=G_0(q,y)+\tilde{G}~,~~~\,
  \tilde{G}:=\sum_{j=1}^{\infty}q^{2j(4j-1)}y^{4j}G_0(q,q^{4j}y)~,$$
  with $|G_0(q,q^{4j}y|\leq |G_0(q,y)|$, see Lemma~\ref{lmG1}.
  This implies (we set $E:=[0,1/2]\times [0,1]$)

  $$\begin{array}{ccl}|G(q,y)|&\geq&\min_{(q,u)\in E}|G_0(q,y)|-
    \sum_{j=1}^{\infty}\max_{(q,u)\in E}(2/3)^{4j}|G_0(q,q^{4j}y)|\\ \\
    &>&g_0^{1/2}-0.85^{1/2}\cdot (16/81)/(1-16/81)~=~0.4402\ldots ~>~0.4402~.
  \end{array}$$

  Part (3) of Lemma~\ref{lmC}. With $G_0$ as above, we set 
  $y:=r(u-\sqrt{1-u^2}i)$,
  $r\in [0,1]$, $u=\cos \phi \in [1/2,\sqrt{2}/2]$, and we find that 
  
  $$\begin{array}{ccl}\partial |G_0|^2/\partial u&=&
    2qr^3(q^8r^4+4q^6r^3u+12q^5r^2u^2-3q^5r^2+q^3r^2+4q^2ru+1)\\ \\
    &\geq&2qr^3(q^8r^4+2q^6r^3+3q^5r^2-3q^5r^2+q^3r^2+2q^2r+1)
    \geq 0\end{array}$$
    (because $2q^5r^2\leq 2q^2r$ and $q^5r^2\leq 1$).
    Equality is possible only for $qr=0$. Next we set $r:=5/7$.
    As in the proof of part (2) we write 
  $G=G_0+\tilde{G}$. Hence
  
    $$\begin{array}{cclcl}|G(q,y)|&\geq&|G_0(q,y)|-
      \sum_{j=1}^{\infty}(5/7)^{4j}|G_0(q,q^{4j}y)|\\ \\
      &\geq&|G_0(q,y)|(~1-(5/7)^4/(1-(5/7)^4)~)&=&|G_0(q,y)|\cdot
      0.6480855856\ldots\end{array}$$
      (we use again the inequality $|G_0(q,q^{4j}y)|\leq |G_0(q,y)|$).
      The modulus $|G_0(q,y)|$
is minimal for $u=1/2$ when it equals

$$(25/5764801)(25q^4+35q^2+49)(625q^8-875q^6+875q^5
-2450q^3-490q^2+1715q+2401)~.$$
One checks numerically that it is an increasing function in~$q$.
Having in mind that $q\geq 0.88$ and setting $G^{\bullet}_0:=|G_0(0.88,(5/7)(1/2-\sqrt{3}i/2)|$, one concludes that 
  
$$|G|\geq G_0^{\bullet}\cdot  0.6480855856\ldots =(0.8459381456\ldots)^{1/2}\cdot
0.6480855856\ldots =0.59607\ldots >0.5960~.$$

  \end{proof}

\begin{proof}[Proof of Lemma~\ref{lmD}]

  The lemma is proved only for the segments with $+$-signs. 
  
  Part (1). We need the following technical result:

\begin{lm}\label{lmFk}
  Consider the product $F_k:=(1+Q/x)(1+qQx)$, $Q=q^k$, $k\in \mathbb{N}$.
  Set $x:=re^{\phi i}$, $\phi \in [0,\pi /2]$.
  Then $\partial |F_k|^2/\partial r\geq 0$ for
  $r\geq q^{-1/2}$, with equality only for $r=q^{-1/2}$.
\end{lm}

\begin{proof}

  We set $u:=\cos \phi$ and hence $\sqrt{1-u^2}=\sin \phi$. Then

  $$\begin{array}{ccll}
    |F_k|^2&=&(2Qru+Q^2+r^2)(Q^2q^2r^2+2qQru+1)/r^3&{\rm and}\\ \\
    \partial |F_k|^2/\partial r&=&2Q(qr^2-1)(Q^2qru+Qqr^2+ru+Q)/r^3~.&
  \end{array}$$
  All factors with all coefficients positive are positive due to $u\geq 0$,
  so $\partial |F_k|^2/\partial r\geq 0$ (resp. $\partial |F_k|^2/\partial r>0$)
  if and only if $r^2\geq 1/q$ (resp. $r^2>1/q$).
\end{proof}

The lemma implies in particular that for $q$ fixed, the quantity
$|\Theta^*|$ increases along the half-line $x=re^{5\pi i/12}$, $r\geq 1.5$.
Hence on the segment $J:=[1.5e^{5\pi i/12},1.7e^{5\pi i/12}]$ one has
$|\Theta^*|<0.021$, see part (1) of Lemma~\ref{lmB} and Lemma~\ref{lmFk}.
Next we subdivide the segment $J$ into subsegments

$$J_1:=[1.5e^{5\pi i/12},1.6e^{5\pi i/12}]~~~\, {\rm and}~~~\,
J_2:=[1.6e^{5\pi i/12},1.7e^{5\pi i/12}]~.$$
Consider the sum
$G_*(q,x):=1/x+q/x^2=(1/x)(1+q/x)$. For $x\in J_1$, one obtains

$$1<|1+q/x|<|1+1/x|\leq |1+(2/3)e^{-5\pi i/12}|=1.3377\ldots <1.3378$$
and $|G_*(q,x)|\leq 1.3378/1.5=0.891866\ldots <0.891867$.
Recall that $q_{\heartsuit}=0.7234985094$. For $q\in [q_{\heartsuit},1)$
  and $x\in J_1$, it is true that

  $$\begin{array}{l}
    1.2175<1.21750\ldots \leq |1+q_{\heartsuit}/x|\leq |1+1/x|~~~\, {\rm and}\\ \\
    |G_*(q,x)|\geq (1/1.6)\cdot 1.2175=0.760937\ldots \geq 0.760937~.
    \end{array}$$
  On the other hand $G=G_*(q,x)+\sum_{j=1}^{\infty}q^{j(2j+1)}G_*(q^{2j+1},x)/x^{2j}$,
  so 

  \begin{equation}\label{equGG}
    \begin{array}{cclcl}
    |G(q,x)|&\geq&|G_*(q,x)|-\sum_{j=1}^{\infty}|G_*(q^{2j+1},x)/x^{2j}|
    &\geq&|G_*(q,x)|-\sum_{j=1}^{\infty}|G_*(q^{2j+1},x)|/(1.5)^{2j}\\ \\
    &\geq&0.760937-0.891867\cdot (1/(1.5^2-1))&=&0.047\ldots ~>~0.021~.
  \end{array}\end{equation}
  As $|G|>|\Theta^*|$, the function $\theta$ has no zeros on~$J_1$. For
  $x\in J_2$ in the same way

  $$|1+q/x|<|1+1/x|\leq |1+(1/1.6)e^{-5\pi i/12}|=1.309255\ldots <1.309256$$
  and $|G_*(q,x)|\leq 1.309256/1.6=0.818285$. For $q\in [q_{\heartsuit},1)$
    and $x\in J_2$, one gets

    $$\begin{array}{l}
    1.1994<1.19939\ldots \leq |1+q_{\heartsuit}/x|\leq |1+1/x|~~~\, {\rm and}\\ \\
    |G_*(q,x)|\geq (1/1.7)\cdot 1.1994=0.705529\ldots \geq 0.705529~.
    \end{array}$$
    By analogy with (\ref{equGG}) one finds that

    $$|G(q,x)|\geq 0.705529-0.818285\cdot (1/(1.6^2-1))=0.18\ldots >0.021$$
    and $\theta$ has no zeros on $J_2$.
\vspace{1mm}

  Part (2). We apply Proposition~\ref{propdecomp} with $n=3$ and $x=re^{\pi i/3}$,
  $0\leq r\leq 1.5$, to express Im$\theta$:
 
  $${\rm Im}(\theta (q,re^{\pi i/3}))={\rm Im}(\theta (q^9,q^{-3}x^3)+
  qx\theta (q^9,x^3)+q^3x^2\theta (q^9,q^3x^3))~.$$
  The quantity $x^3=-r^3$ is real negative and each of the functions
  $\theta (q^9,q^sx^3)$, $s=-3$, $0$ or $3$, is real-valued. Clearly
  Im$(\theta (q^9,-q^{-3}r^3))=0$. Therefore, as
  $\sin (\pi /3)=\sin (2\pi /3)=\sqrt{3}/2$,  
 
  $${\rm Im}(\theta (q,re^{\pi i/3}))=(\sqrt{3}/2)(qr\theta (q^9,-r^3)
  +q^3r^2\theta (q^9,-q^3r^3))~.$$
  We make use of a result in \cite{Kosepar} which says that $\theta$
  has no real zeros $\geq -6$. Given that $\theta (q,0)=1$, this means that
 $\theta (q,x)>0$ for $x\in [-6,0]$. Hence $\theta (q^9, -r^3)>0$ for
 $r\leq 6^{1/3}$ and $\theta (q^9,-q^3r^3)>0$ for $qr\leq 6^{1/3}$. As
 $1.5<6^{1/3}$, this means that Im$(\theta (q,re^{\pi i/3}))>0$ for $r\leq 1.5$,
  so $\theta$ has no zeros in the given segments.
  \vspace{1mm}

  Part (3). We apply Proposition~\ref{propdecomp} with $n=4$
  and $x=re^{\pi i/4}$, $0\leq r\leq 1.4$:

  $${\rm Im}(\theta (q,re^{\pi i/4}))={\rm Im}(\theta (q^{16},q^{-6}x^4)+
  qx\theta (q^{16},q^{-2}x^4)+
  q^3x^2\theta (q^{16},q^2x^4)+q^6x^3\theta (q^{16},q^6x^4))~.$$
  Clearly $x^4=-r^4$, Im$(\theta (q^{16},q^{-6}x^4))=0$ and

  $${\rm Im}(\theta (q,re^{\pi i/4}))=(\sqrt{2}/2)(qr\theta (q^{16},-q^{-2}r^4)+
  q^6r^3\theta (q^{16},-q^6r^4))+q^3r^2\theta (q^{16},-q^2r^4)~.$$
  The functions $\theta (q^{16},-q^{\ell}r^4)$, $\ell =-2$, $2$ or $6$,
  are positive-valued for $-q^{\ell}r^4\in (-6,0)$. For $r\in (0,1.4]$ and
 $q\in [0.82,1]$, one has indeed $-q^{\ell}r^4\in (-6,0)$, so
 Im$(\theta (q,x))\neq 0$.

 Suppose that $q\in [0.7,0.82]$, with $0.7<q_{\heartsuit}$. One has $-q^2r^4\in (-6,0)$ and
 $-q^6r^4\in (-6,0)$, so $\theta (q^{16},-q^2r^4)>0$,
 $\theta (q^{16},-q^6r^4)>0$ and it remains to consider $\theta (q^{16},-q^{-2}r^4)$. As 
 $q^{16}\in (0.003,0.042)$ and as by \cite[Proposition~7]{KoBSM1},
 $\theta (q,x)>0$ for $x\in [-1/q,0]$, one obtains that
 $\theta (q^{16},-q^{-2}r^4)>0$ for $-q^{-2}r^4>-0.042^{-1}=-23.8\ldots$.
 However for $(q,r)\in [0.7,0.82]\times (0,1.4]$, one has
$-q^{-2}r^4\geq -7.84$, so $(\theta (q^{16},-q^{-2}r^4))>0$ and
Im$(\theta (q,x))>0$.
\vspace{1mm}

Part (4) results from \cite[Theorem~1]{KoMatStud}.
\end{proof}

\section{Proof of Theorem~\protect\ref{tmimprovebis}
  \protect\label{sectmimprovebis}}

For $q\in (0,\tilde{q}_1]$, the function $\theta$ has no complex zeros
and all its real zeros are $<-5$, see \cite[Proposition~2]{KoMatStud}.
Therefore we prove Theorem~\ref{tmimprovebis}
only for $q\in (\tilde{q}_1,1)$.

\subsection{Proof for $q\in (\tilde{q}_1,0.44]$}

We deduce the proof from the following proposition:

\begin{prop}
  (1) For $q\in (\tilde{q}_1,0.44]$, there is no zero of $\theta$
    on the segments $[3e^{\pm 3\pi i/4},4.3e^{\pm 3\pi i/4}]$.
 
 (2) For $q\in (\tilde{q}_1,0.44]$, there is no zero of $\theta$ on the arc 
 $\{ |x|=4.3$, {\rm arg}$x\in [3\pi /4,5\pi /4]\}$.
\end{prop}

Indeed, the zeros of $\theta$ depend continuously on $q$.
There are no zeros of $\theta$ on the rest of the border
of the exclusion zone with Re$x\leq 0$, see part~(3)
of Theorem~\ref{tmoldresults}. Hence no zero enters this zone
for $q\in (\tilde{q}_1,0.44]$. The proposition is proved only for Im$x\geq 0$.

\begin{proof}[Proof of part (1)]
  We set $x=r(-\sqrt{2}+i\sqrt{2})/2$, $r\in [3,4.3]$,
  and $\theta :=\theta_4+\theta^{\dagger}$,
  with $\theta_4:=\sum_{j=0}^4q^{j(j+1)/2}x^j$. Then
  Im$(\theta_4)=(qr)\theta^{\diamond}$, where
  $\theta^{\diamond}:=1/\sqrt{2}-q^2r+q^5r^2/\sqrt{2}$. Notice that 

  $$\partial \theta^{\diamond}/\partial r=-(-2q^3r+\sqrt{2})q^2/\sqrt{2}~~~\,
  {\rm and}~~~\, \partial \theta^{\diamond}/\partial q=
  -(-5q^3r+2\sqrt{2})qr/\sqrt{2}~.$$
  Both derivatives are negative for
  $(q,r)\in [\tilde{q}_1,0.44]\times [3,4.3]$. Hence 

  $${\rm Im}(\theta_4)\geq \tilde{q}_1\cdot 3\cdot
  \theta^{\diamond}(0.44,4.3)=0.0837\ldots ~.$$
  At the same time $|\theta^{\dagger}|\leq \theta^{\spadesuit}:=
  \sum_{j=5}^{\infty}0.44^{j(j+1)/2}4.3^j=0.0068\ldots <0.0837$.
  Thus Im$(\theta )\neq 0$, so $\theta \neq 0$.
\end{proof}

\begin{proof}[Proof of part (2)]
 We set $x:=4.3\cdot(u+\sqrt{1-u^2}i)$, $u\in [-1,-\sqrt{2}/2]$. Then 
 
 $${\rm Im}(\theta_4)=\sqrt{1-u^2}P~,~~~\, {\rm where}~~~\,
 P|_{u=-1}=4.3q+238.521q^6-36.98q^3-1367.5204q^{10}~.$$
 It is easy to check numerically that $P|_{u=-1}$
 is decreasing in $q$, with $P_{u=-1,q=0.44}=0.1007\ldots$. 
 
 One finds directly that 
 
 $$\partial P/\partial u=q^3R~,~~~\,
 R:=36.98+636.056q^3u+8205.1224q^7u^2-1367.5204q^7~.$$
 The inequality $R>0$ is a corollary of
 $R|_{u=-1}=36.98-636.056q^3+6837.602q^7$ (which is positive-valued
 on $[\tilde{q}_1,0.44]$) and $\partial R/\partial u=q^3(636.056+16410.2448q^4u)$
 which is also positive-valued. Thus $R>0$ and $P\geq 0.1007$. 
 
 For $u\in [-0.9949,-\sqrt{2}/2]$, one has $\sqrt{1-u^2}>0.1$
 and Im$(\theta_4)\geq 0.01007$ while $|\theta^{\dagger}|\leq \theta^{\spadesuit}$.
 Hence Im$(\theta )\neq 0$ and $\theta \neq 0$.
 So suppose that $u\in [-1,-0.9949]$. We minorize the function
 
 $$\begin{array}{ccl}{\rm Re}(\theta_4)&=&
   1-18.49q^3+341.8801q^{10}+4.3qu-238.521q^6u\\ \\
   &&+36.98q^3u^2-2735.0408q^{10}u^2+318.028q^6u^3+2735.0408q^{10}u^4\end{array}$$
   by substituting $-1$ for $u$ in the fourth, seventh and eighth monomial
   and $-0.9949$ in the fifth, sixth and ninth. The polynomial thus obtained
   takes only values $>0.1>\theta^{\spadesuit}$ in $[\tilde{q}_1,0.44]$.
   This proves part~(2).
\end{proof}

\subsection{Beginning of the proof for $q\in [0.44,1)$ --
    the function $\Theta^*$\protect\label{subsecbeginning}}

We need a technical lemma:

\begin{lm}\label{lmfg}
  (1) Set $f:=|(1+qr(u+\sqrt{1-u^2}i))(1+q(u-\sqrt{1-u^2}i)/r)|$, where
  $u:=\cos \phi$, $\phi \in [0,\pi ]$ and
  $r>0$. For fixed $q$ and $r$, the quantity $f$
  is an increasing function in $u$ on $[-1,1]$.
  \vspace{1mm}

  (2) For $r\in [3,4.3]$ and $\phi \in [3\pi /4,\pi ]$, the quantity
  $g:=|1+(u-\sqrt{1-u^2}i)/r|$ is maximal if and only if
  $\phi =3\pi /4$, $r=4.3$ when it equals $g^{***}:=0.8515846485\ldots$.
  \end{lm}

\begin{proof}
  Part (1) of the lemma follows from

  $$\begin{array}{ccl}f^2&=&1+q^2r^2+2q^3ru+q^2/r^2+q^4+2qu/r+
    4q^2u^2+2qru+2q^3u/r~~~\, {\rm and}\\ \\ 
    \partial f^2/\partial u&=&
    2q((q^2r^2+2qru+1)+(r^2+2qru+q^2))/r\geq 0~,\end{array}$$
  with equality only for $r=q=1$, $u=-1$ (quadratic polynomials
  with non-positive discriminants).

  Part (2) results from $g^2=1+2u/r+1/r^2$ which for fixed $r$, is maximal for
  $u=-\sqrt{2}/2$; for $u=-\sqrt{2}/2$, $g^2$ is maximal for $r=4.3$.  
\end{proof}

With the help of the lemma we find a majoration of the quantity $|\Theta^*|$.
By Lemma~\ref{lmfg}, it suffices to consider the case $u=-\sqrt{2}/2$.
We study the product

$$Y(q,r):=(1-q)(1+qre^{3\pi i/4})(1+qe^{-3\pi i/4}/r)~.$$
Clearly
$\Theta^*(q,re^{3\pi i/4})=(1+e^{-3\pi i/4}/r)\prod_{j=1}^{\infty}Y(q^j,r)$. 

\begin{lm}
  It is true that $|Y(q,r)|^2=(1-q)^2(q^2r^2-\sqrt{2}qr+1)M(q,r)$, where
  $M(q,r)=(r^2-\sqrt{2}qr+q^2)/r^2$. 
  \end{lm}

The lemma is to be checked straightforwardly.

\begin{lm}\label{lm013}
  For $(q,r)\in [0,1]\times [3,+\infty )$, one has
  $\partial M/\partial r>0$ and $\partial M/\partial q<0$. Hence if
    $q\in [\tilde{q}_1,1)$ and $r\in [3,4.3]$, then the product
      $\prod_{j=1}^{\infty}M(q^j,r)$ is maximal for
      $(q,r)=(\tilde{q}_1,4.3)$ when it equals
    $0.8632817769\ldots$. 
\end{lm}

\begin{proof}
  A direct computation shows that
  $\partial M/\partial r=(\sqrt{2}qr-2q^2)/r^3>0$ and
  $\partial M/\partial q=(2q-\sqrt{2}r)/r^2<0$.  
  \end{proof}

\begin{lm}\label{lmchi}
  (1) Suppose that $3\leq r_1<r_2\leq 4.3$. Set $\chi_r(q):=q^2r^2-\sqrt{2}qr+1$.
  Then for $q\in [0,1]$, the graphs of the functions $\chi_{r_j}(q)$, $j=1$, $2$,
  intersect only for $q=0$ and $q=q^{**}:=\sqrt{2}/(r_1+r_2)$. For $q\in (0,q^{**})$
  (resp. $q\in (q^{**},1]$), one has $\chi_{r_1}>\chi_{r_2}$ (resp. $\chi_{r_1}<\chi_{r_2}$).  
\vspace{1mm}

(2) For $q\in [0,1]$ and $r\in [3,4.3]$, one has $0\leq (1-q)^2\chi_r(q)\leq 1$,
with equalities only for respectively $q=1$ and $q=0$.

(3) For $q\in (0,1)$, it is true that
$\tilde{\chi}(q):=\max_{r\in [3,4.3]}\chi_r(q)=\max (\chi_3(q),\chi_{4.3}(q))$.
\end{lm}

 \begin{proof}
   The first statement of part (1) is immediate. The other statements
   follow from $\chi_{r_2}(1)=r_2(r_2-\sqrt{2})+1>\chi_{r_1}(1)=r_1(r_1-\sqrt{2})+1$
   and $\chi_{r_1}'(0)=-\sqrt{2}r_1>\chi_{r_2}'(0)=-\sqrt{2}r_2$. 

   Part (2). First of all, $\chi_r(q)=(qr-\sqrt{2}/2)^2+1/2>0$.
   The quadratic function 
$$f_r^{\diamondsuit}(q):=(qr-\sqrt{2}/2)(1-q)$$
takes its maximum at $r_0:=(2r+\sqrt{2})/4r\in (1/2,1)$, with 
$$f_r^{\diamondsuit}(r_0)=(-2r+\sqrt{2})^2/16r\in [0.43,0.76]~.$$
The function is non-negative on $[\sqrt{2}/2r,1]$ and negative
on $[0,\sqrt{2}/2r)$, with $f_r^{\diamondsuit}(0)=-\sqrt{2}/2=-0.707\ldots$.
  For $q\in (0,0.3)$, $|qr-\sqrt{2}/2|\leq
  \max (0.3\cdot 4.3-\sqrt{2}/2,\sqrt{2}/2)=\sqrt{2}/2$, so 

  $$(f_r^{\diamondsuit}(q))^2\leq 1/2~~~\, {\rm and}~~~\,
  (1-q)^2\chi_r(q)=(f_r^{\diamondsuit}(q))^2+(1/2)(1-q)^2\leq (1/2)(1+(1-q)^2)\leq 1~.$$
  For $q\in [0.3,1]$, $(f_r^{\diamondsuit}(q))^2\leq 0.76^2$ and
  $(1-q)^2\chi_r(q)\leq 0.76^2+(1/2)\cdot 0.7^2<1$.

Part (3) results from part~(1).
 \end{proof}
 
 We set $\kappa (q):=(1-q)^2\tilde{\chi}(q)$,
 $\kappa^{\omega}(q):=\prod_{j=1}^{\omega}\kappa (q^j)$, where
 $\omega \in \mathbb{N}^*\cup \infty$.  From part (2) of Lemma~\ref{lmchi}
 one deduces that $\kappa^{\infty}(q)<\kappa^9(q)$. It is clear that
 for $r\in [3,4.3]$, $|Y(q^j,r)|^2\leq (1-q^j)^2\tilde{\chi}(q^j)M(q^j,r)$,
 therefore 

 $$|\Theta^*(q,re^{3\pi i/4})|^2\leq |1+e^{-3\pi /4}/r|^2\cdot
 \kappa^{\infty}(q)\cdot \prod_{j=1}^{\infty}M(q^j,r)<\kappa^9(q)~,$$
 see part (2) of Lemma~\ref{lmfg} and Lemma~\ref{lm013}.
 One checks numerically that for $q\in [0.55,1)$, the function $\kappa^9$
   takes only values $<0.031$. Hence
   $|\Theta^*(q,re^{3\pi i/4})|<0.031^{0.5}=0.176\ldots$. In the case
   $q\in [0.44,0.55)$ the estimation of $|\Theta^*|$
     is carried out in Subsection~\ref{subseccompletion}.

     \subsection{The proof for $q\in [0.44,1)$ -- the
         function $G$\protect\label{subsecG}}

       In this subsection we minorize the quantity $|G(q,x)|$.
       In the proofs of the lemmas lowercase indices indicate
       the corresponding partial derivations. For instance,
       $S_r$ and $S_{r^4}$ stand for $\partial S/\partial r$ and
  $\partial^4S/\partial r^4$. 
  
\begin{lm}\label{lmLr}
  Set $x:=r(u+\sqrt{1-u^2}i)$, $u:=\cos \phi$,
  $\phi \in [3\pi /4,\pi ]$, $r\geq 3$, $H:=1/x+q/x^2+q^3/x^3+q^6/x^4$
  and $L:=|H|^2$.
  Then $\partial L/\partial r<0$.
\end{lm}

\begin{proof}
Using computer algebra one
  finds that $L_r=-2S/r^9$, where

  $$\begin{array}{cclccl}
    S&=&r^6+3r^5uq+K_4r^4+K_3r^3+K_2r^2+7q^9ru+4q^{12}~,&
    K_4&=&8u^2q^3-4q^3+2q^2~,\\ \\ K_3&=&20u^3q^6-15uq^6+5q^4u~,&
    K_2&=&12q^7u^2-6q^7+3q^6~.\end{array}$$
  We prove that $S>0$ from which the lemma follows. We find first

  $$S_{r^3}|_{r=3}=6(20u^3q^6-15uq^6+5q^4u+96q^3u^2-48q^3+24q^2+270uq+540)~.$$
  It is positive-valued due to $|270qu|\leq 270$, $|20u^3q^6|\leq 20$,
  $|5q^4u|\leq 5$, $|-48q^3|\leq 48$ and $270+20+5+48<540$. Then we consider

  $$S_{r^4}=96(2u^2-1)q^3+360r(uq+r^2)+48q^2$$
  which is positive, because $2u^2\geq 1$ and $|uq|<r^2$. Thus $S_{r^3}>0$.
  Next,

  $$S_{r^2}|_{r=3}=360q^6u^3+(27q^7+864q^3)u^2+(-270q^6+90q^4+1620q)u-12q^7+
  6q^6-432q^3+216q^2+2430$$
  which is also positive-valued due to $27q^7u^2-12q^7\geq 0$,
  $|360q^6u^3|\leq 864q^3u^2$ and
  $$|(90q^4+1620q)u|+|-432q^3|\leq 90+1620+432=2142<2430~.$$
  This together with
  $S_{r^3}>0$ implies $S_{r^2}>0$. Then we study

  $$\begin{array}{ccl}S_r|_{r=3}&=&
    540q^6u^3+(72q^7+864q^3)u^2+(7q^9-405q^6+135q^4+1215q)u\\ \\
    &&-36q^7+18q^6-432q^3+216q^2+1458~.\end{array}$$
  Summing up the inequalities

  $$\begin{array}{cclccl}
    540(q^6u^3+q^3u^2)&\geq&0~,&72q^7u^2-36q^7&\geq&
    0~,\\ \\ 324q^3u^2-432q^3+216q^2+101&>&0~,&
    (7q^9+135q^4+1215q)u+1357&\geq&0~,\\ \\ &{\rm and}&&-405q^6u+18q^6&\geq&0~,
  \end{array}$$
  one gets $S_r|_{r=3}>0$. Combined with $S_{r^2}>0$, this gives $S_r>0$. Next we
  deal with

  \begin{equation}\label{equSr3}\begin{array}{ccl}S|_{r=3}&=&
    540q^6u^3+(108q^7+648q^3)u^2+(21q^9-405q^6+135q^4+729q)u\\ \\ &&
    +4q^{12}-54q^7+27q^6-324q^3+162q^2+729~.
  \end{array}\end{equation}
  Suppose first that $u\in [-0.95,-\sqrt{2}/2]$. Some of the terms
  participating in (\ref{equSr3}) will be split in sums as follows:
  $540q^6u^3=432q^6u^3+108q^6u^3$, $729=693+36$ and
  $648q^3u^2=319q^3u^2+108q^3u^2+221q^3u^2$. One has

  $$|729qu|\leq 729\cdot 0.95<693~,~~~\,
  156\sqrt{2}=220.6\ldots <221~~~\, {\rm and}~~~\, 
  |(21q^9+135q^4)u|\leq 221q^3u^2~~~(B)~.$$
  The inequality $S|_{r=3}>0$ results from $(B)$ and 

  $$\begin{array}{cclccl}
    729qu+693&>&0~,&(108u^2-54)q^7&\geq&0~,\\ \\
    (432u^3-405u+27)q^6&\geq&0~,&108(q^6u^3+q^3u^2)&\geq&0~,\\ \\
    -36q^3+36&\geq&0~,&319q^3u^2+4q^{12}-288q^3+162q^2&\geq&0~.
    \end{array}$$
  Now we suppose that $u\in [-1,-0.95]$. Then one obtains a minoration of
  $S|_{r=3}$ by substituting in (\ref{equSr3}) $-1$ for $u$ and $u^3$,
  and $0.95^2$ for $u^2$. The values of the thus obtained polynomial
  in $q$ are $>206$ for $q\in (0,1)$. Thus $S|_{r=3}>0$
  which together with $S_r>0$ implies $S>0$. The lemma is proved.

\end{proof}

\begin{lm}\label{lmLu}
  For $q\in (0,1)$, $u\in [-1,-\sqrt{2}/2]$ and $r\geq 3$,
  one has $\partial L/\partial u>0$.
\end{lm}

\begin{proof}
 One finds that 
 
 $$\partial L/\partial u=2qL^*/r^7~,~~~\,
 L^*:=12q^5r^2u^2+(4q^6r+4q^2r^3)u+q^8-3q^5r^2+q^3r^2+r^4~.$$
 As $6q^5r^2u^2-3q^5r^2\geq 0$, to show that $L^*>0$ and hence
 $\partial L/\partial u>0$, it suffices to prove that 
 
 $$L^{**}:=6q^5r^2u^2+(4q^6r+4q^2r^3)u+q^8+q^3r^2+r^4>0~.$$
 We notice first that 
 
 $$\partial L^{**}/\partial r=12q^5ru^2+(4q^6+12q^2r^2)u+2q^3r+4r^3>0~,$$
 because $4q^6u+2q^3r\geq 0$ and $12q^2r^2u+4r^3\geq 0$.
 So to prove the lemma it suffices to show that 
 
 $$L^{**}|_{r=3}=54q^5u^2+(12q^6+108q^2)u+q^8+9q^3+81>0~.$$
 Suppose first that $q^2\leq 3/4$. Then $12q^6u+54q^5u^2\geq 0$ and
 $108q^2u+81\geq 81u+81\geq 0$, so $L^{**}|_{r=3}>0$. Suppose now that
 $q^2\geq 3/4$. We observe that $6q^5u^2+12q^6u+9q^3\geq 0$ and that
 $48q^5u^2\geq 48\cdot(9/16)\cdot qu^2=27qu^2$, so 
 
 $$L^{**}|_{r=3}\geq 27qu^2+108q^2u+q^8+81>27qu^2
 +108\sqrt{q}u+81=27(\sqrt{q}u+1)(\sqrt{q}u+3)\geq 0$$
 from which the lemma follows.
\end{proof}

We observe that $G(q,x)=\sum_{k=0}^{\infty}(q^{8k^2-2k}/x^{4k})H(q,x/q^{4k})$, so

$$|G|\geq |H(q,x)|-\tilde{H}~,~~~\,
\tilde{H}:=\sum_{k=1}^{\infty}|(q^{8k^2-2k}/x^{4k})H(q,x/q^{4k})|~. $$
We minorize $|H(q,x)|$ by setting $x=-4.3$, see
Lemmas~\ref{lmLr} and \ref{lmLu}. Thus $|H(q,-4.3)|=-H(q,-4.3)$ and

\begin{equation}\label{equHq}
  d|H(q,-4.3)|/dq=-1/4.3^2+3q^2/4.3^3-6q^5/4.3^4<-1/4.3^2+3/4.3^3<0~.
\end{equation}
In what follows we consider the situation when $r\geq 3$ and $q\in (0,q_{**}]$, 
$0<q_{**}\leq 1$, so 

$$|H(q,-4.3)|\geq |H(q_{**},-4.3)|=:g(q_{**})~.$$

     Next we majorize $\tilde{H}$. We set $A:=|q^{8k^2-2k}/x^{4k}|$, so
     $A\leq (q_{**})^{8k^2-2k}/3^{4k}$.
     To majorize $|H(q,x/q^{4k})|$ we make use of the following lemma:

\begin{lm}\label{lm|H|}
For $r\geq 3$ fixed, one has $H^{\ast}:=d|H(q,r(-\sqrt{2}+\sqrt{2}i)/2)|^2/dq<0$.
\end{lm}

\begin{proof}
  We have $H^{\ast}=(-r^{-8})S$, where 

  $$S:=-12q^{11}+9\sqrt{2}rq^8-6\sqrt{2}r^3q^5+4\sqrt{2}r^3q^3-6r^2q^5
  +\sqrt{2}r^5-2r^4q~.$$
  Obviously $-12q^{11}+9\sqrt{2}rq^8>0$,
  $-6\sqrt{2}r^3q^5+4\sqrt{2}r^3q^3>-2\sqrt{2}r^3q^5$ and

  $$-2\sqrt{2}r^3q^5-6r^2q^5
  +\sqrt{2}r^5-2r^4q>-2\sqrt{2}r^3-6r^2+\sqrt{2}r^5-2r^4=:
  r^2\tilde{\delta}(r)~.$$
  As $\tilde{\delta}^{(j)}(3)>0$, $j=0$, $1$, $2$ and $3$, we have $\tilde{\delta}(r)>0$, $S>0$ and
  $H^{\ast}<0$ for $r\geq 3$.
\end{proof}
Suppose that $r\geq 3$ and $q\in (0,q_{**}]$. By Lemmas~\ref{lmLr} and
   \ref{lm|H|},

    $$ |H(q,x/q^{4k})|\leq |H(q,3/(q_{**})^{4k})|\leq |H(0,3/(q_{**})^{4k})|=
     (q_{**})^{4k}/3~.$$
   Thus $\tilde{H}\leq \sum_{j=1}^{\infty}
   |((q_{**})^{8k^2-2k}/3^{4k})((q_{**})^{4k})/3)|=
   \sum_{j=1}^{\infty}|(q_{**})^{8k^2+2k}/3^{4k+1}|=:\mu (q_{**})$ and

   \begin{equation}\label{equgmu}
     |G|\geq g(q_{**})-\mu (q_{**})~.
     \end{equation}
\vspace{1mm}

{\em Conclusion.} Suppose that $(r,u)\in [3,4.3]\times [-1,-\sqrt{2}/2]$.
We use inequality (\ref{equgmu}) to minorize $|G|$ for different values of
$q_{**}$. Lemmas~\ref{lmLr}, \ref{lmLu} and inequality (\ref{equHq})
imply that 

\begin{enumerate}
\item 
  for $q\in [0,0.55]$,
  the minimum of $|H|$ is attained for $q=0.55$
  which equals $g_1:=g(0.55)=0.2048239481\ldots$; one has
  $\mu_1:=\mu (0.55)=0.000010\ldots$ and $|G|\geq 0.2048$; 

\item for $q\in [0.55,1]$,
  the minimum of $|H|$ is attained for $q=1$ and equals 
$g_2:=g(1)=0.1881273580\ldots$; one has $\mu_2:=\mu (1)=0.00416666\ldots$ and 
  $|G|\geq 0.1839$;
  
\item for $q\in (0,0.47]$, the minimum of $|H|$ is attained for
 $q=0.47$ which equals $g_3:=g(0.47)=0.2084132995\ldots$; one has
$\mu_3:=\mu (0.47)=2.16\ldots \times 10^{-6}$ and $|G|\geq 0.2084$; 

\item for $q\in (0,0.445]$, the minimum of $|H|$ is attained for
 $q=0.445$ which equals $g_4:=g(0.445)=0.2095767068\ldots$; one has
$\mu_4:=\mu (0.445)=1.25\ldots \times 10^{-6}$ and $|G|\geq 0.2095$.

\end{enumerate}

\subsection{Completion of the proof for
  $q\in [0.44,1)$\protect\label{subseccompletion}}

  Suppose first that $q\in [0.55,1)$. At the end of
    Subsection~\ref{subsecbeginning} we obtained the estimation
    $|\Theta^*|<0.031^{0.5}=0.176\ldots$. Having in mind that
    $|G|\geq 0.1839$ (see case~2. the end of Subsection~\ref{subsecG}),
    we conclude that $|\theta |\geq |G|-|\Theta^*|>0.006>0$,
    so $\theta \neq 0$.

    Suppose next that $q\in [0.47,0.55]$. We remind that
    the quantities $M$ and $\kappa^9(q)$ were introduced
    in Subsection~\ref{subsecbeginning}. The quantity $M$
    is maximal for $q=0.47$, with 
    $$M=0.7475993990\ldots =:M^{\diamond}~~~\, {\rm and}~~~\,
    (M^{\diamond})^{1/2}=0.8646383053\ldots ~.$$
    One finds that $\kappa^9(q)\leq \kappa^9(0.47)=0.06328020076\ldots$
    and $(\kappa^9(0.47))^{1/2}=0.2515555620\ldots =:\kappa^{\diamond}$.
    Thus $|\Theta^*(q,4.3e^{3\pi i/4})|\leq \kappa^{\diamond}\cdot
    (M^{\diamond})^{1/2}\cdot g^{***}=0.1848\ldots<0.2048\leq |G|$,
    see Lemma~\ref{lmfg} and case~1.
    at the end of Subsection~\ref{subsecG}. Hence again $\theta \neq 0$.
\vspace{1mm}

Suppose that $q\in [0.445,0.47]$. Then $M$ is maximal
for $q=0.445$, with 
$$M=0.7686911306\ldots =:M^{\bullet}~~~\, {\rm and}~~~\,
(M^{\bullet})^{1/2}=0.8767503240\ldots ~.$$
It is true that $\kappa^9(q)\leq \kappa^9(0.445)=0.07530778828\ldots$
and $(\kappa^9(0.445))^{1/2}=0.2744226453\ldots =:\kappa^{\bullet}$. Thus 
$|\Theta^*(q,4.3e^{3\pi i/4})|\leq \kappa^{\bullet}\cdot
(M^{\bullet})^{1/2}\cdot q^{***}=0.2045\ldots <0.2084\leq |G|$,
see case~3. at the end of Subsection~\ref{subsecG}.
As above we conclude that $\theta \neq 0$.
\vspace{1mm}

Suppose that $q\in [0.44,0.445]$. Then $M$ is maximal
for $q=0.44$, with 
$$M=0.7727506376\ldots =:M^{\circ}~~~\, {\rm and}~~~\,
(M^{\circ})^{1/2}=0.8790623627\ldots ~.$$
Next, $\kappa^9(q)\leq \kappa^9(0.44)=0.07786987357\ldots$
and $(\kappa^9(0.44))^{1/2}=0.2790517400\ldots =:\kappa^{\circ}$. Hence 
$|\Theta^*(q,4.3e^{3\pi i/4})|\leq \kappa^{\circ}\cdot
(M^{\circ})^{1/2}\cdot q^{***}=0.2085\ldots$. At the same time $|G|$
is minoirized by $0.2095$ (see case~4. at the end of
Subsection~\ref{subsecG}), so once again $\theta \neq 0$.

\end{document}